\pdfoutput=1
\documentclass[final]{siamart0516}
\usepackage{graphicx}
\usepackage{amssymb}
\usepackage{amsmath}
\usepackage{subfig}
\usepackage{multirow}

\newcommand{\nrlev}{\ensuremath{L}}  
\newcommand{\nrsam}{\ensuremath{M}}  
\newcommand{\nrBZp}{\ensuremath{K}}  
\newcommand{\supsz}{\ensuremath{N}}  
\newcommand{\Pvac}{\ensuremath{p}}   
\newcommand{\energy}{\ensuremath{\epsilon}}
\newcommand{\smooth}{\ensuremath{\delta}}
\newcommand{\IDoS}{\ensuremath{I}}
\newcommand{\DoS}{\ensuremath{\rho}}
\newcommand{\e}[1]{\ensuremath{\times 10^{#1}}}
\newcommand{\E}[1]{{\ensuremath{\mathbb{E}}\mspace{-2mu}\left[#1\right]}}
\def\VAR{\mathrm{Var}}  
\newcommand{\tol}{\textrm{TOL}}      
\newcommand{\rset}{\mathsf{R}}
\newcommand{\zset}{\mathsf{Z}}
\def\MoS{\mathrm{MoS}}
\def\AMC{\mathcal{A}_{\mathrm{MC}}}
\def\AMLMC{\mathcal{A}_{\mathrm{MLMC}}}

\headers{MLMC for randomly perturbed materials}{P. Plech\'a\v{c}, E. von Schwerin}

\title{Multi-level {M}onte {C}arlo acceleration of computations on
  multi-layer materials with random defects\thanks{%
  Version: \today.
\funding{This work was funded by the  U.S. DOD-ARO Grant Award W911NF-14-1-0247.}
}
}

\author{
  Petr Plech\'{a}\v{c}\thanks{Department of Mathematical Sciences, University of Delaware, Newark, DE, 19716 
  {(\email{plechac@udel.edu})}},    
  \and
  Erik von Schwerin\thanks{Department of Mathematical Sciences, University of Delaware, Newark, DE, 19716 
  {(\email{schwerin@udel.edu})}}.
}

\begin{document}

\maketitle

\begin{abstract}
  We propose a Multi-level Monte Carlo technique to accelerate
  Monte Carlo sampling for approximation of properties of 
  materials with random defects.
  The computational efficiency is investigated on test problems
  given by tight-binding models of a single layer of graphene or of 
  $\MoS_2$ where the integrated electron density of states per
  unit area is taken as a representative quantity of interest.
  For the chosen test problems the multi-level Monte Carlo estimators
  significantly reduce the computational time of standard Monte Carlo
  estimators to obtain a given accuracy. 
\end{abstract}

\section{Introduction}
\label{sec:introduction}

The aim of this study is to develop non-intrusive numerical techniques for
approximating properties of layered heterostructures 
with impurities in random locations. 

The goal is to apply these techniques on models developed and used for
layered heterostructures 
such as tight-binding models for transition-metal
dichalcogenides (TMDC). 
The numerical techniques are
not in themselves restricted to tight-binding models, but can be
combined with more computationally intensive and accurate models 
when such are called for. 
For the purpose of testing and calibrating the algorithms we use two 
simple tight-binding models of materials with honeycomb lattices. The
first is of graphene, where individual atoms at 
random locations are ``removed'' from the tight-binding model without
changing the positions of the surrounding atoms. This example
can be viewed as a rough approximation of a graphene sheet where hydrogen
atoms are attached to randomly distributed carbon atoms in the sheet without
mechanically deforming the sheet. We also use a tight-binding model of
a single layer of the TMDC $\MoS_2$; in this model
$\mathrm{S}$ atoms are similarly removed.

Characteristically we wish to compute material properties which, in
the unperturbed case of a periodically repeating fundamental cell, can
be obtained by computing the band structure of the material over the
first Brillouin zone associated with the fundamental
cell.
Here we test the approach on computations of the integrated electronic
density of states per unit area of the material, which can be computed
from the band structure and is a common quantity of interest in such
computations. 
This is interesting on its own, and also serves as a test case for the
more computationally demanding problem of computing
the electric conductivity by the Kubo-Greenwood
formula. 
This tensor depends both on the energies of the band
structure and on the corresponding eigenstates.

We assume that the random perturbations result in an ergodic random field.
Random perturbations of the studied material break the periodicity,
which is used when computing the band structure and quantities
depending upon it. A common approach 
in this case is to 
extend the fundamental cell of the unperturbed material along the
primitive lattice vectors. 
In the test case this means extending the fundamental cell of the
honeycomb lattice by some integer factors $\supsz_1$ and $\supsz_2$ along its
primitive lattice vectors. Random perturbations are 
introduced in this ``super cell'' of the fundamental cell, which is
then periodically extended to 
cover the whole plane. The band structure can now be computed, but at
a much higher cost, increasing with the size of the super
cell. Finally, in theory, the size of the super cell is allowed to go
to infinity to obtain the limit of random perturbations without
periodicity.
In the remainder of this paper we will let $\supsz_1=\supsz_2=\supsz$.

The discrete random perturbations in our test examples only allow a
finite number of outcomes for each finite super cell.
Indeed, if the super cell is small enough it is efficient to compute
the quantity of interest for all possible combinations of
perturbations, which with the known probability of each outcome gives
a complete description of the random quantity of interest.
This exhaustive approach can not be
extended to large super cells where our true interest lies in the case
of a randomly perturbed material; neither can it be applied where the
perturbations vary continuously. 
We are left with approximate methods.
For low densities of impurities in particular,
homogenisation or perturbation around the solution of the problem
without impurities can be used. 
Here we will instead study Monte Carlo
methods to estimate the expected value of the quantity of interest,
which is a random variable for a finite size super cell. The main
advantage of Monte Carlo type methods is in their simplicity; they are
non-intrusive methods in the sense that they approximate the expected
value of the desired quantity by the average over several independent
outcomes of the random perturbations, and each outcome can be computed
using any existing code capable of taking the perturbed configuration as
an input.

Our goal is to use so called Multilevel Monte Carlo methods to
reduce the computational cost of standard Monte Carlo sampling while
retaining the same accuracy.
The key point here is to systematically generate control variates to
an expensive, accurate, numerical approximation of a random sample. 
With a suitable choice of control variates fewer samples on the most
expensive and most accurate approximation level are needed and the
total computational cost to reach a given accuracy can be reduced.
In Section~\ref{sec:MC_MLMC} we will describe Monte Carlo and
Multilevel Monte Carlo estimators for the test problem and discuss the
numerical complexity assuming simplified models of the approximation
properties and computational cost of the underlying computational
method.
In Section~\ref{sec:tb} we describe our tight-binding test problems and
explain how to generate control variates for Multilevel Monte Carlo in
this context.
Finally, in Section~\ref{sec:Numerics} we show numerical experiments
which illustrate the efficiency of the multilevel approach on the
given test problems.

\section{Monte Carlo and MultiLevel Monte Carlo}
\label{sec:MC_MLMC}

By Monte Carlo methods here we simply mean methods where the expected
value of a random variable is approximated by the sample average over
several, usually independent, outcomes of the random variable. 
In the present context this means that we generate a number of
outcomes of the random perturbations of the materials model on the
super cell and then compute the quantity of interest for each outcome
individually by separate calls to the underlying computational physics
code. 
In this spirit we want to restrict ourselves to Monte Carlo techniques
that do not strongly depend on the particular qualities of our test
problem; for example we do not in this discussion optimize our methods
given test problem by utilizing the 
fact that only a finite number of perturbations are possible for each
finite super cell. 

\subsection{Monte Carlo complexity}
\label{sec:MC}

The quantity of interest in the test problem applied in 
Section~\ref{sec:tb}, which is an integrated density of states,
is a deterministic quantity in the infinite volume limit, $|V|\sim
\supsz^d\to\infty$; 
that is the variance goes to zero as the
size of the super cell goes to infinity. Does this mean that we should
use only one sample of the random perturbations in the material? 

We can answer the above question by estimating the
rate at which the variance of our 
quantity of interest goes to zero as the super cell size, $\supsz$,
increases, and compare this to the rate at which the expected value of the
quantity converges and the rate at which the computational work grows. 
Let $Q$ be the exact value, in this case deterministic, of the
quantity we wish to approximate, let $Q_\supsz$ be the random variable of
the same quantity computed on a finite super cell of size $\supsz$ with
random perturbations, and let $\E{Q_\supsz}$ and $\VAR(Q_\supsz)$ denote the
expected value and the variance of $Q_\supsz$, respectively.
Assume the following models for the above quantities:
\begin{subequations}
  \label{eq:compl_model}
  \begin{align}
    \label{eq:Q_bias}
    \text{the finite $\supsz$ bias, } && Q-\E{Q_\supsz}&\propto \supsz^{-W},\\
    \label{eq:Q_variance}
    \text{the variance, } && \VAR(Q_\supsz)&\propto \supsz^{-S},\\
    \label{eq:Q_cost}
    \text{the cost per sample, }&& \mathtt{cost} & \propto \supsz^{C},
  \end{align}  
\end{subequations}
for some positive constants $W$, $S$, $C$.
Assume, for now, that the number of samples, $\nrsam\to\infty$, and
approximate the expected value $\E{Q_\supsz}$ by the estimator
\begin{align}
  \label{eq:MC}
  \AMC(\nrsam) & = \frac{1}{\nrsam}\sum_{m=1}^\nrsam Q_\supsz(m),
\end{align}
where $Q_\supsz(m)$ denotes the $m$:th independent sample of $Q_\supsz$. Then 
by the Central Limit Theorem we can justify approximating the
suitably rescaled \emph{statistical error} of our Monte Carlo
estimator by the Standard Normal random variable, $N(0,1)$, which allows us to
state the following error constraints. To make the total error in our
quantity of interest approximately $\tol$ with high probability, we
require that the bias is approximately $(1-\Theta)\tol$ for some
$\Theta\in(0,1)$ and the variance of our Monte Carlo estimator is
approximately $\frac{1}{C_\alpha}(\Theta\tol)^2$ where the confidence
parameter $C_\alpha$ is chosen for a Standard Normal random
variable. That is 
\begin{align*}
  \supsz^{-W} & \approx (1-\Theta)\tol,\\
  \frac{1}{\nrsam}\supsz^{-S} & \approx \frac{1}{C_\alpha}(\Theta\tol)^2.
\end{align*}
Minimizing the total work, proportional to $\supsz^{C}$, with respect to
$\Theta$ while satisfying the two constraints leads to the simple
expression for the splitting parameter 
\begin{align*}
  0<\Theta&=\frac{1}{1+\frac{C-S}{W}}<1,
\end{align*}
provided that the cost of generating samples grow faster than the
variance of the sampled random variables decrease,
i.e. $C>S$. Furthermore, the optimal number of samples becomes
$\nrsam\propto\tol^{-(2-S/W)}$ which, as long as $S<2W$, goes to infinity as
$\tol\to 0$. 
With the work per sample assumed to be $\supsz^C$ and with $\supsz^{-W} \approx
(1-\Theta)\tol$ the total work for a Monte Carlo method is then
approximately proportional to 
\begin{align}
  \label{eq:Work_SLMC}
  \mathtt{Work}_\mathtt{MC}(\tol) & \propto \tol^{-(2+\frac{C-S}{W})}.
\end{align}

A method using a fixed number of samples must take
$\supsz\propto\tol^{-2/S}$, assuming that $S<2W$, giving the
asymptotic complexity
\begin{align}
  \label{eq:Work_FS}
  \mathtt{Work}_\mathtt{FS}(\tol) & \propto \tol^{-\frac{2C}{S}}.
\end{align}
Thus, the Monte Carlo complexity~\eqref{eq:Work_SLMC} is an
improvement as long as $C>S$. 

Qualitatively the above argument tells us that for small error
tolerance it is more computationally efficient to use several samples
on a smaller super cell than to use a larger super cell with only one
sample of the random perturbations. 
For quantitative predictions on the optimal choice we may use a
sequence of increasing super cell sizes to empirically estimate the
parameters in the models for how the bias and variance decays with $\supsz$
and how the work grows with $\supsz$. From these estimates we can decide
how to optimally choose the number of samples versus the size of the
super cell.

\subsection{Multilevel Monte Carlo as an acceleration of 
  standard Monte Carlo}
\label{sec:MLMC}

Assume that the models~\eqref{eq:compl_model} hold approximately for
large enough $\supsz$ and that parameters, $W$, $S$, and $C$, have been
empirically or theoretically estimated and found to be such that it is
more efficient to use Monte Carlo sampling than one single sample on a very
large super cell. 
In this situation we want to use Monte Carlo methods to approximate
the expected value of a quantity which in turn has a bias due to a
method parameter; in this case we assume most importantly by the size
given test problem of the super cell, $\supsz$. 
Over the past decade so called Multilevel Monte Carlo (MLMC) method
has become an 
increasingly popular systematic technique for accelerating such 
Monte Carlo methods. They can be traced
back to Heinrich et al.~\cite{heinrich98,hs99} where they were
introduced for parametric integration, and were independently proposed
by Giles~\cite{giles08} in a form closer to the one in this
paper. Following~\cite{giles08} the methods have typically been
applied to problems where each sample of a 
standard Monte Carlo sample is obtained by the solution of a
discretization based numerical approximation to a stochastic
differential equation or a partial differential equation with random data. 
This technique
depends on the possibility of using cheaper approximations of
the quantity to be evaluated for each random sample as control
variates for more accurate approximations; see~\cite{GilesAcNum}. 
For example, in a discretization based numerical method characterized
by a mesh size, $h$, with known convergence as $h\to 0$, a solution
using a larger step size $2h$ can be used as a control variate to a
solution using a step size $h$ which have been chosen to make the bias
sufficiently small. A good use of control variates means that fewer
samples on the accurate, most expensive, scale can be used, while
samples on less accurate and less costly scales are introduced to
compensate.

In the present context the artificially finite super cell size
introduces a bias which only vanishes in the limit as $\supsz\to\infty$. We
also assume that among the parameters in the numerical approximation
$\supsz$ dominates the computational cost as our tolerated error $\tol\to
0$. 
It is then natural to consider using approximate values of our
quantity of interest based on smaller super cell sizes as control
variates to the more accurate approximations computed on large super
cells. 
Assume, for now, that for $\supsz_\ell=c2^\ell$, with $c,\ell\in\zset_+$, in
addition to the approximate quantity of interest $Q_\ell$ on super cell
size $\supsz$ we can construct control variates $Q_\ell^{CV}$ such that 
\begin{subequations}
  \label{eq:compl_model2}
  \begin{align}
    \label{eq:CV_unbiased}
    \E{Q_\ell^{CV}} & = \E{Q_{\ell-1}},\\
    \label{eq:CV_variance}
    \VAR(Q_\ell-Q_\ell^{CV}) & \propto \supsz^{-D},
  \end{align}
\end{subequations}
for some $D>S$, 
and the cost of sampling the control variate is small compared to
sampling the original quantity of interest; at most a constant fraction
smaller than one say, so that~\eqref{eq:Q_cost} holds for generating
the pair $(Q_\ell,Q_\ell^{CV})$. 
Following the standard MLMC approach the estimator~\eqref{eq:MC} is
now replaced by
\begin{align}
  \label{eq:MLMC}
  \AMLMC & = \frac{1}{\nrsam_1}\sum_{m=1}^{\nrsam_1}Q_{1}(\omega_{1,m}) + 
  \sum_{\ell=2}^\nrlev\frac{1}{\nrsam_\ell}
  \sum_{m=1}^{\nrsam_\ell}\left(Q_{\ell}(\omega_{\ell,m})-Q_{\ell}^{CV}(\omega_{\ell,m})\right),
\end{align}
where $\supsz_\ell=c2^\ell$ for $\ell=1,2,\dots,\nrlev$, and $\nrsam_\ell$ denotes the
positive integer number of samples used on size $\supsz_\ell$; by
$\omega_{\ell,m}$ we denote the $m$:th independent identically
distributed outcome of the random impurities on a super cell of size
$\supsz_\ell$. 
Note that while we assume independence between all terms of the sums
in~\eqref{eq:MLMC}, the difference
$Q_{\ell}(\omega_{\ell,m})-Q_{\ell}^{CV}(\omega_{\ell,m})$
is computed using the same outcome of the random perturbation but two
different approximations of $Q$.

Taking the expected value, the sum over $\ell$ in the
definition~\eqref{eq:MLMC} telescopes by assumption~\eqref{eq:CV_unbiased}
so that $\AMLMC$ is an unbiased estimator of
$\E{Q_{\nrlev}}$. Furthermore, by independence of the outcomes
$\omega_{\ell,m}$,
\begin{align*}
  Var{\left(\AMLMC\right)} & = 
  \frac{1}{\nrsam_1}\VAR(Q_{1}) +
  \sum_{\ell=2}^\nrlev\frac{1}{\nrsam_\ell}\VAR(Q_{\ell}-Q_{\ell}^{CV}),
\end{align*}
where the variances are assumed approximated by~\eqref{eq:Q_variance}
and~\eqref{eq:CV_variance}. 
Similarly to the standard Monte Carlo case we require that the sum of
the bias and the statistical error of the estimator sum up to a
specified error tolerance, $\tol$. 
Denote by $W_\ell$ the work, as 
modeled by~\eqref{eq:Q_cost}, of computing one sample
on level $\ell$, that is $Q_1$, for $\ell=1$, or
$Q_{\ell}-Q_{\ell}^{CV}$, for $\ell=2,\dots,\nrlev$. 
Also let $V_\ell$ denote the corresponding variances predicted by
models~\eqref{eq:Q_variance}, for $\ell=1$,
and~\eqref{eq:CV_variance}, for $\ell=2,\dots,\nrlev$. 
A straightforward minimization of the computational work model with
respect to the number of samples on each level leads to
\begin{align}
  \label{eq:opt_samples}
  \nrsam_\ell & = \left( \frac{C_\alpha}{ \theta\tol} \right)^2
  \sqrt{\frac{V_\ell}{W_\ell}} \sum_{k=1}^\nrlev \sqrt{W_k V_k},
  && \text{for $\ell=1,\dots,\nrlev$}
\end{align}
in terms of general work estimates, $\{W_\ell\}_{\ell=1}^\nrlev$, and
variance estimates, $\{V_\ell\}_{\ell=1}^\nrlev$; see for
example~\cite{Haji_Ali_opt_hier}. Here the number of levels, $\nrlev$,
depends on $\tol$ through the constraint on the finite $\supsz$ bias.

Further minimizing the predicted work of generating
$\AMLMC$ 
with respect to the splitting between bias and statistical error, the
model of the computational work becomes
\begin{align}
  \label{eq:Work_MLMC}
  \mathtt{Work}_\mathtt{MLMC}(\tol) & \propto \tol^{-(2+\frac{C-D}{W})}.
\end{align}
This improves on the computational work of a standard Monte Carlo
method as long as $D>S$, that is as long as
$\VAR(Q_{\ell}-Q_{\ell}^{CV})$ decays at a higher rate in $\supsz_\ell$
than $\VAR(Q_{\ell})$. The applicability of MLMC techniques depends
on finding control variates satisfying this condition. We will
describe how to generate such control variates in Section~\ref{sec:CV}.

\section{Tight-binding model with random defects}
\label{sec:tb}

In our test problems the target is to compute the
integrated density of states in tight-binding models of a single-layer
material with honeycomb lattices. 
The first example is a simple nearest neighbor
tight-binding model of graphene, which provides us with a well
controlled, and comparatively inexpensive, 
test setting where we can study the algorithms before turning to new
materials.
The second example is a tight-binding model of one layer of
$\MoS_2$. 

\subsection{Materials model without defects}
\label{sec:model_unpert}

In a tight-binding model of a periodically repeating material, we take
a given numbering of the atoms in the fundamental cell of the periodic
material and identify periodic images of the atoms. 
Using values for hopping and on-site energies obtained for example by
parameter fitting to more accurate density functional theory results
we construct a Hamiltonian matrix, $H(k)$, and an overlap matrix,
$S(k)$, leading to a generalized eigenvalue problem
\begin{align}
  \label{eq:gen_evp}
  H(k)u & = \energy S(k)u.
\end{align}
Our quantities of interest will depend on the solutions
to~\eqref{eq:gen_evp} for each point $k$ in the Brillouin zone.

\paragraph{A tight-binding model for graphene}

Here we use a 
nearest neighbor tight-binding model of a single-layer graphene
sheet from Chapter 2.3.1, ``$\pi$~Bands of Two-Dimensional Graphite'',
in~\cite{carbon_nanotubes}. 

In this tight-binding model, including only the $\pi$ energy bands,
the generalized eigenvalue problem~\eqref{eq:gen_evp} is defined by
\begin{subequations}
  \label{eq:tightbinding}
  \begin{align}
    \label{eq:tb_Hdiag}
    H_{mm}(k) & = \energy_{2p} \\
    \label{eq:tb_Hcross}
    H_{mn}(k) & = t_{mn}\exp{(ik\cdot R_{mn})} \\
    \label{eq:tb_Sdiag}
    S_{mm}(k) & = 1 \\
    \label{eq:tb_Scross}
    S_{mn}(k) & = s_{mn}\exp{(ik\cdot R_{mn})}
  \end{align}  
\end{subequations}
where $R_{mn}$ is the
vector from atom position $m$ to $n$ in the honeycomb lattice. In the
nearest neighbor interactions the parameters $t_{mn}$ and $s_{mn}$
are 0 unless atoms $m$ and $n$ are nearest neighbors and
$t_{mn}=\langle\phi_m|\mathcal{H}|\phi_n\rangle=t$ and
$s_{mn}=\langle\phi_m,\phi_n\rangle=s$, independent of
$m$ and $n$, otherwise. The numerical values were taken
from~\cite{carbon_nanotubes} to be
\begin{align*}
  \energy_{2p}&=0\,\mathrm{eV}& t&=-3.033\,\mathrm{eV} 
  & s &=0.129\,\mathrm{eV}
\end{align*}
which gives the Fermi level $\energy_F=0\,\mathrm{eV}$.

The fundamental cell of the honeycomb lattice of the graphene sheet has
two atoms, call them $A$ and $B$, so that $H(k)$ and $S(k)$ are 2-by-2
matrices where by the periodic structure the only non-diagonal elements
$H_{AB}(k)=H_{BA}(k)^\star$ are obtained by
summing~\eqref{eq:tb_Hcross} over the three nearest neighbor
directions; similarly $S_{AB}(k)=S_{BA}(k)^\star$ is obtained
from~\eqref{eq:tb_Scross}. 

\paragraph{A tight-binding model of $\MoS_2$}

In an ideal single layer $\MoS_2$, the projection of the atom
positions on the plane forms a honeycomb lattice, just as for
graphene. This time the two types of lattice positions, $A$ and $B$,
are occupied by an $\mathrm{Mo}$-atom and a pair of
$\mathrm{S}$-atoms, separated in the direction perpendicular to the
plane of the  $\MoS_2$ layer; see Figure~\ref{fig:CV}.

In this example we take the tight-binding model of a monolayer TMDC
material from Section~IV, equations~(4)--(10),
in~\cite{PhysRevB.92.205108}, and the parameters for $\MoS_2$
in Table~VII of the same paper. 
This model includes 11 bands and interactions up to selected
third-neighbor couplings which together define the Hamiltonian
$H(k)$; the overlap matrix $S(k)$ is the identity matrix. 

\subsection{Materials model with defects}
\label{sec:model_pert}

We now consider the case when individual atom locations in the
infinite sheet of the honeycomb lattice are ``removed'' from the
tight-binding model. 
In the graphene case, 
we view this as a rough approximation to hydrogen
atoms attaching to the corresponding carbon atoms and thus changing
the electron interactions without mechanically deforming the sheet.
Still in the graphene case, the locations of the removed atom locations
are random, and it is assumed that each atom location is removed, with
probability $\Pvac$, $0<\Pvac<1$, independently of all other locations.
A vacancy is modeled by removing all rows and columns corresponding to
interactions involving this lattice site from the tight-binding
Hamiltonian matrix, $H(k)$, and overlap matrix, $S(k)$. 

\begin{figure} \centering
  \includegraphics[width=0.3\textwidth]{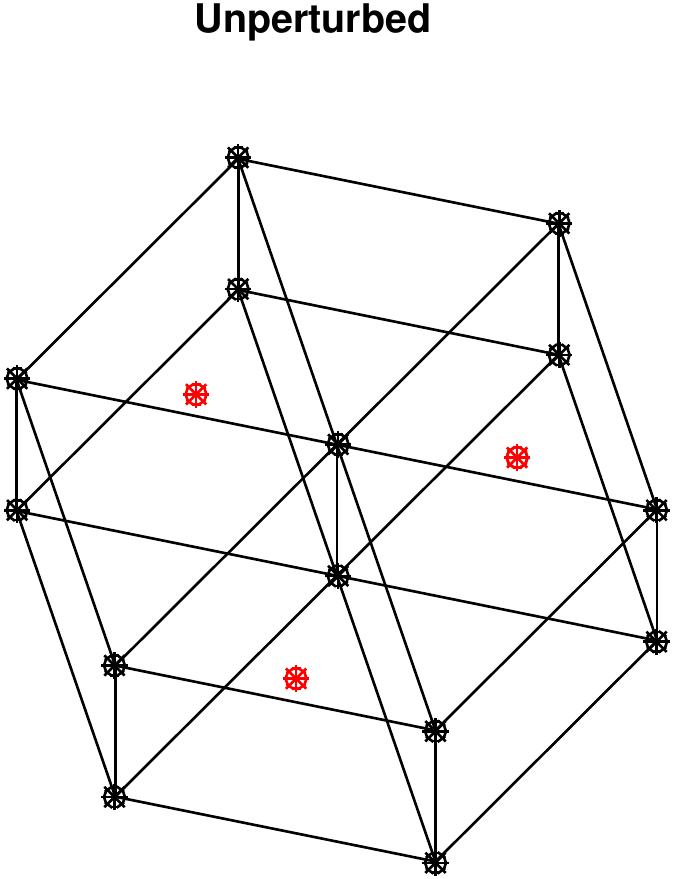}
  \includegraphics[width=0.3\textwidth]{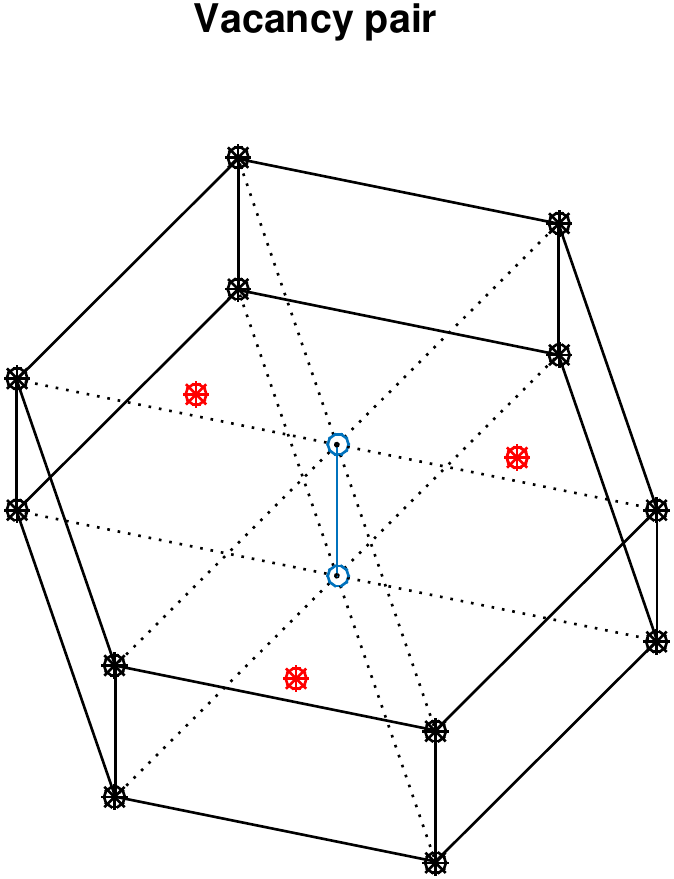}
  \caption{$\mathrm{MoS_2}$: (Left) Unperturbed MoS2 (Right) Perturbed
    by vacancy pair.}
  \label{fig:Viz_MoS2}
\end{figure}

In a simplified test of perturbations of the $\MoS_2$ layer,
in order to keep the example similar to the graphene model, we let
the permutations \emph{remove pairs} of $\mathrm{S}$ atoms
located at randomly sites, instead of individual S
atoms; see Figure~\ref{fig:Viz_MoS2}. Any such pair of
$\mathrm{S}$ atoms is removed with probability $\Pvac$ 
independently of all other pairs. No Mo atoms are removed. 
The numerical tests include three different probabilities
$\Pvac=0.025$, 0.05, and 0.1.
Such a vacancy pair is modeled by removing from the tight-binding
Hamiltonian matrix, $H(k)$, all rows and columns
corresponding to Wannier orbitals involving this particular pair of
$\mathrm{S}$ atoms. 

Ideally the atom locations should be chosen independently of each
other on the entire infinite sheet, but as described above, this is
approximated by taking a large super cell where the impurities are
distributed randomly; this super cell and its impurities are then
repeated periodically to create an infinite sheet. 
We seek the limit as the size of the super cell goes to infinity, and
commit an approximation error by keeping it finite.

\subsection{Control Variates for an MLMC Approach}
\label{sec:CV}

\begin{figure}
  \centering
  \includegraphics[width=0.8\textwidth]{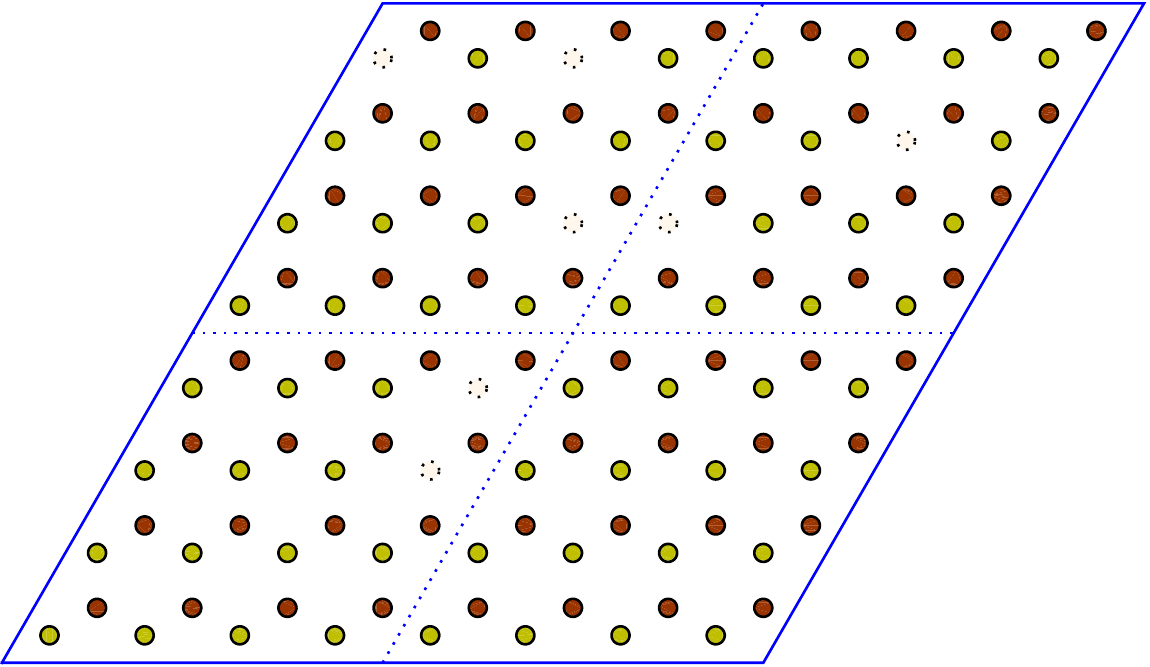}
  \includegraphics[width=0.4\textwidth]{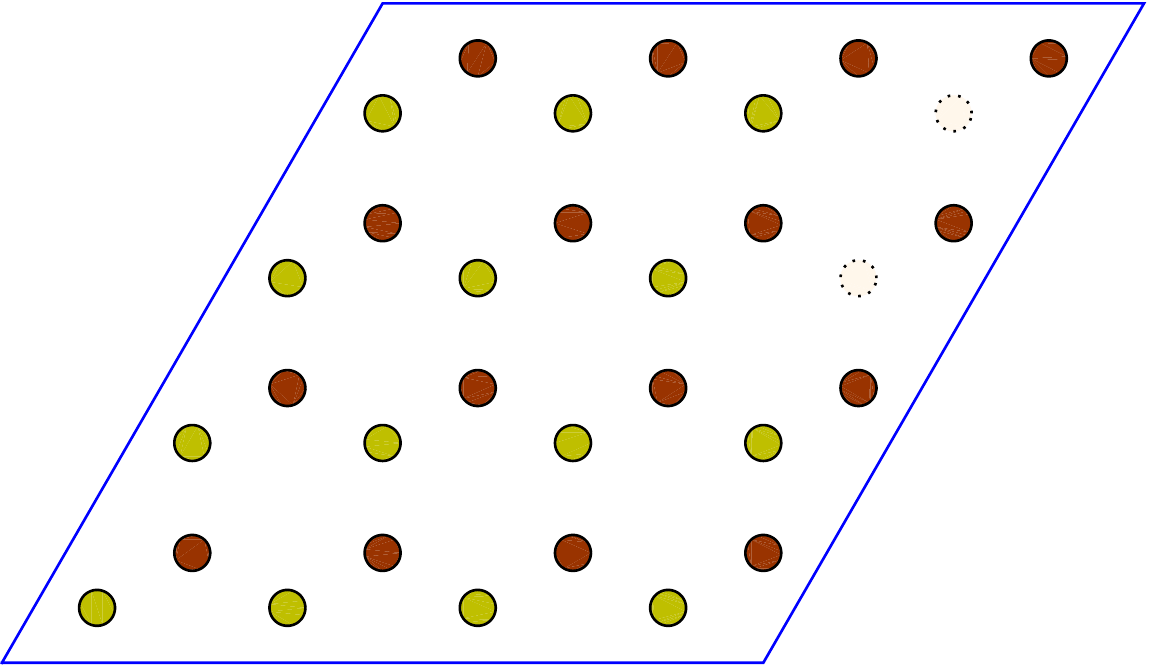}
  \includegraphics[width=0.4\textwidth]{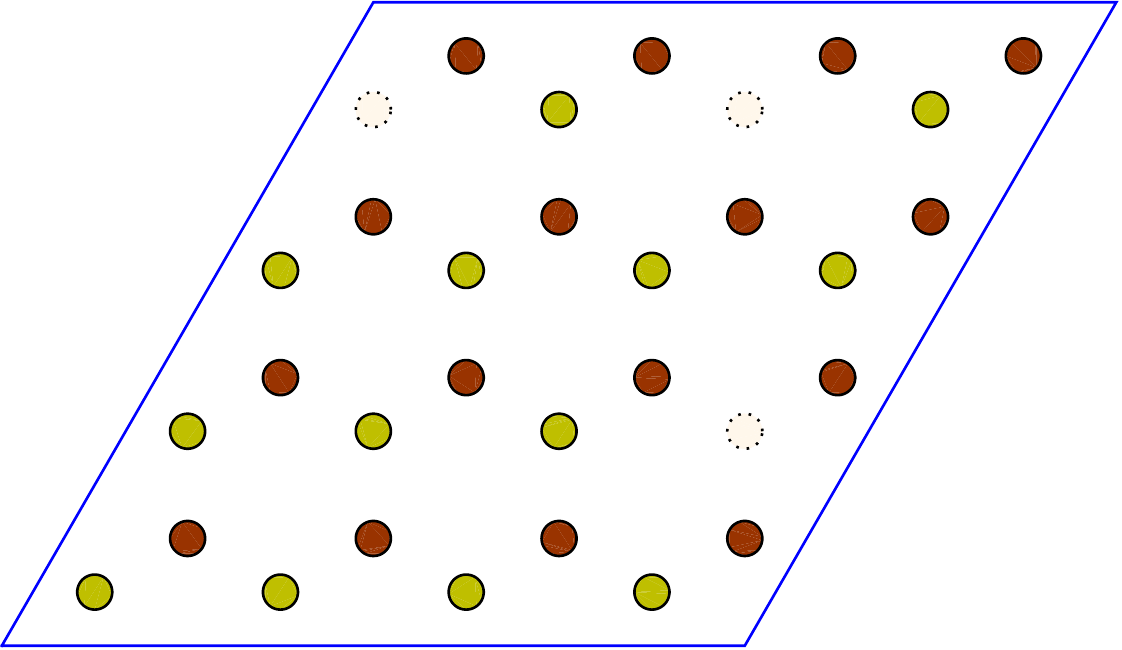}\\
  \includegraphics[width=0.4\textwidth]{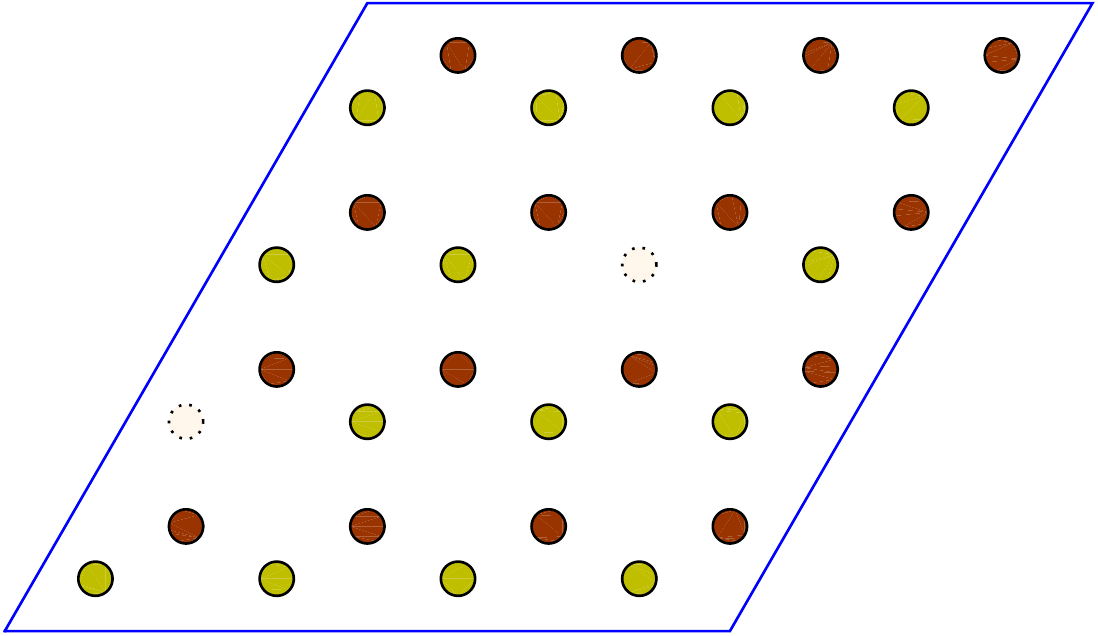}
  \includegraphics[width=0.4\textwidth]{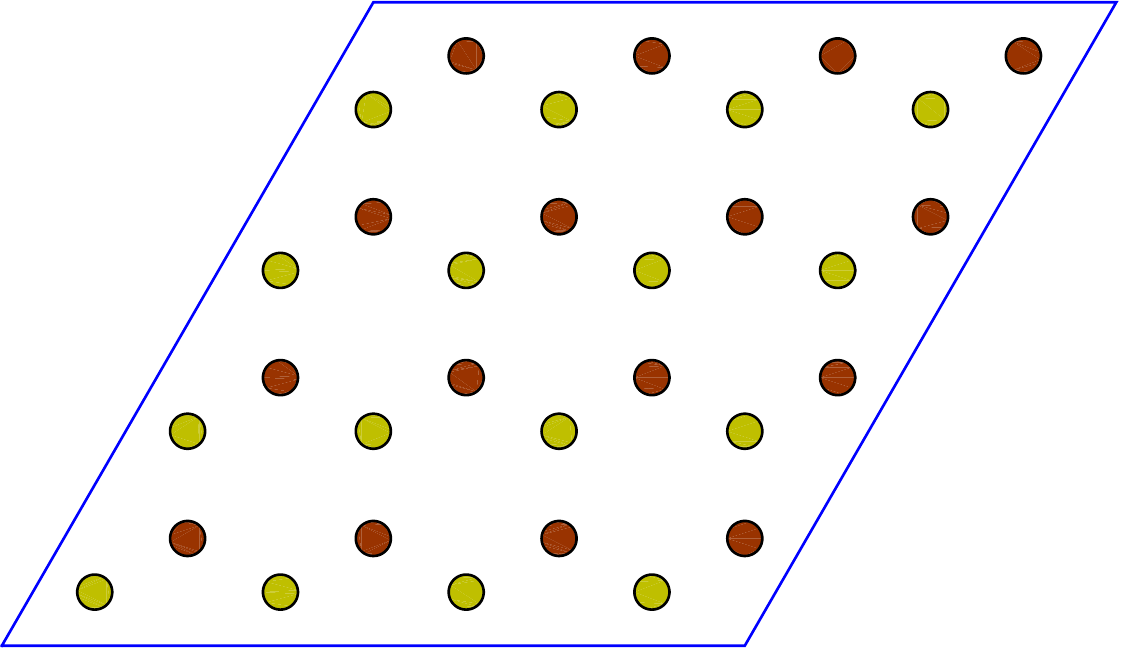}
  \caption{Control variate idea. In the 8-by-8 super cell on the top
    the brown and yellow circles illustrate atom sites of type ``A''
    and ``B'', respectively. Seven circles have been
    removed indicating that the corresponding sites have been replaced
    by vacancies.
    This is one random outcome of the impurities
    on a rather small super cell.
    The larger super cell has been divided into four corners 
    which
    all inherit the impurities of the corresponding part of the larger
    super cell. These four smaller super cells are themselves extended
    periodically to the entire plane; the quantity of interest is
    computed on all four of them, and the arithmetic mean is used as
    a control variate for the quantity computed on the larger super
    cell.}
  \label{fig:CV}
\end{figure}

The MLMC approach to accelerate the standard Monte Carlo sampling
introduced in Section~\ref{sec:MC_MLMC} rests on the possibility to
automatically generate control variates for the random variable whose
expected value we wish to approximate. The control variates must be 
cheaper to sample than the target random variable while
still being strongly correlated to the target.
In our randomly perturbed tight-binding model the dominating factor in
the computational cost of generating one sample is the size of the
finite super cell, $\supsz$. It is thus natural to try control variates on
smaller super cells which, for any given outcome of random impurities,
resemble the larger super cell. Assume for example that $\supsz$ is
divisible by 2. We can divide a large super cell into four parts
where each part retains the impurities of the larger super cell as
illustrated in Figure~\ref{fig:CV} and then extend each part
periodically to an infinite sheet. The quantity 
of interest computed on each one of the four parts will be correlated
to that computed on the larger super cell, and we can take the
arithmetic mean of the four parts as our control variate.

More generally, let $\mathcal{F_\ell}$ denote the $\ell$:th supercell
in the MLMC hierarchy,
$P(\ell)$ denote the number of atom sites in $\mathcal{F_\ell}$, and
$X=(x_1,\dots,x_{P(\ell)})$ be the coordinates of the $P(\ell)$ atom
sites.
We represent a partition of $\mathcal{F_\ell}$ into $R$ subdomains by
the function $\Phi_\ell:\mathcal{F_\ell}\to\{1,\dots,R\}$.
We then define the control variate 
\begin{align*}
  Q_{\ell}^{CV}(\omega;\mathcal{F_\ell}) & = 
  \frac{1}{R}\sum_{r=1}^RQ_{\ell-1}\left(\omega;\Phi_\ell^{-1}(r)\right),
\end{align*}
where $\omega$ denotes a the outcome of the random perturbation on
level $\ell$ and $Q_{\ell-1}\left(\omega;\Phi_\ell^{-1}(r)\right)$
denotes the quantity of interest computed on the subproblem restricted
to $\Phi_\ell^{-1}(r)$. 
We require that the partition is chosen so that
$Q_{\ell-1}\left(\dot;\Phi_\ell^{-1}(r)\right)$ are i.i.d. random
variables for independent outcomes of the random perturbations to
guarantee that condition~\eqref{eq:CV_unbiased} is satisfied.
In the specific case of the tight-binding models in
Section~\ref{sec:model_pert}, this restricted subproblem involves
solving generalized eigenvalue problems~\eqref{eq:gen_evp} with
matrices $H(k)$ and $S(k)$ satisfying the periodicity condition on the
new subdomains.

This systematic way of generating control variates in a multilevel
approach can be naturally extended to other geometries, for example 
an infinite nano ribbon. The random
impurities could then model either impurities along the edge following some
given distribution or again atoms binding to the surface of the ribbon
in random locations. The requirement~\eqref{eq:CV_unbiased} will be
satisfied as long as the super cell in this quasi 1D problem is
divided along the direction of the ribbon.

\section{Numerical Tests}
\label{sec:Numerics}

Here we empirically investigate whether the proposed control variates
satisfy the conditions under which MLMC improves on the computational
complexity of standard Monte Carlo sampling.

\subsection{Quantities of Interest}
\label{sec:QoI}

The physical quantity to approximate from our computational
model in the test case is the integrated electronic density of
states of the material. 
For a periodic material, before we let the
artificial finite size of the super cell go to infinity, this
property depends on the bandstructure computed over the first
Brillouin zone.

\subsection{Numerical approximation of bandstructure}
\label{sec:bandstructure}


The first Brillouin zone associated with the fundamental cell of the
honeycomb lattice is a regular hexagon. 
For the unperturbed material, it is by symmetry sufficient to consider
a rhombus which constitutes one third of the Brillouin zone. This
rhombus is here uniformly divided into $\nrBZp_1$ by $\nrBZp_2$ rhombi, with 
discretization points, $k_{mn}$, in the corners of the rhombi. For
each $k_{mn}$ the generalized eigenvalue problem~\eqref{eq:gen_evp} is
solved numerically using Matlab's full eigenvalue solver ``\texttt{eig}''.

\begin{figure}
  \centering
  \includegraphics[width=0.45\textwidth]{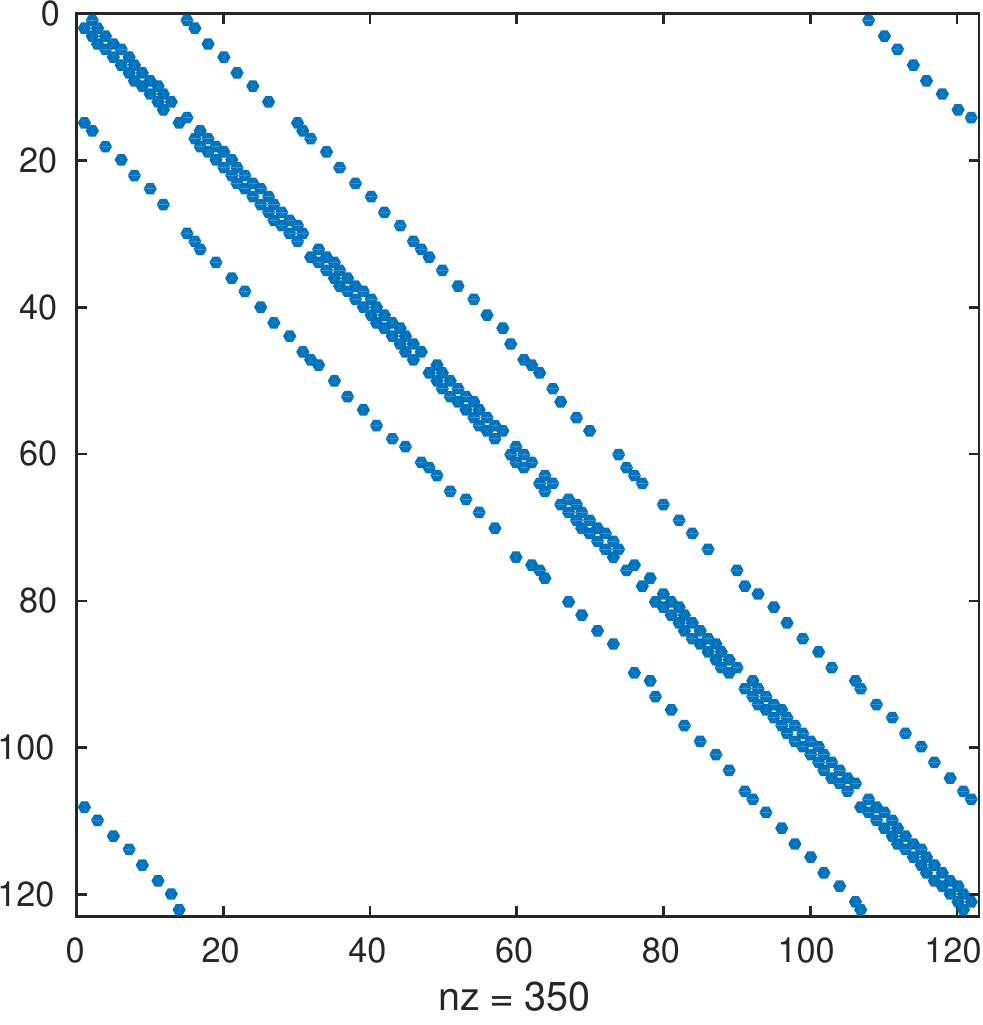}
  \includegraphics[width=0.45\textwidth]{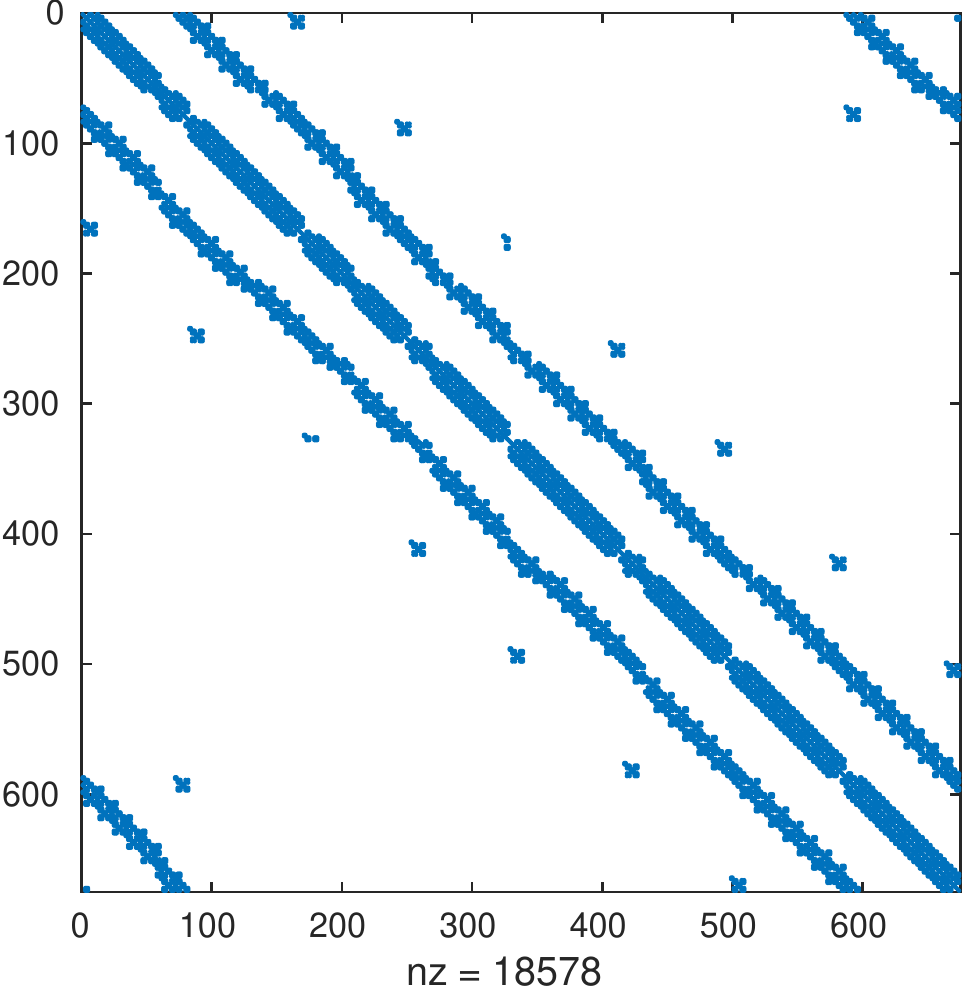}
  \caption{Sparsity structure of an outcome of the matrix $H$ with
    $\supsz=8$ for the graphene (left) and $\MoS_2$ (right) models.}
  \label{fig:sparsity}
\end{figure}

Note that for a nearest neighbor tight-binding model the matrices of
the generalized eigenvalue problem are sparse; see
Figure~\ref{fig:sparsity} for examples with $\supsz=8$. 
As $\supsz$ grows larger one must take advantage of the sparsity in the
eigenvalue computations. However, more complex tight-binding models
will be less sparse, and in more accurate density functional theory
computations the corresponding problems become non-linear and very
much more complex to solve.

For a super cell where the fundamental cell has been extended by an
integer factor $\supsz$ along both primitive lattice vectors, the first
Brillouin zone is still a regular hexagon, rescaled by the factor
$1/\supsz$. Perturbations in random atom locations in the periodically
repeating super cell break the symmetry which allowed us to compute the
bandstructure on a third of the Brillouin zone. Hence the
bandstructure is computed on three rhombi which combined make up the
Brillouin zone.
In all the numerical examples we used $\nrBZp_1=\nrBZp_2=\nrBZp$,
where in turn the resolution in the Brillouin zone was held constant
as $\supsz$ increased; that is $\supsz\nrBZp=\mathtt{constant}$. In
the graphene example $\nrBZp=64/\supsz$ and in the $\MoS_2$
example $\nrBZp=128/\supsz$.

\subsection{Numerical approximation of the integrated density of
  states} 
\label{sec:IDoS}

The quantity of interest in the present test is the expected value of
the integrated density 
of states. The electronic density of states per unit area of the
two-dimensional material, $\DoS(\energy)$ at energy
$\energy$, is defined as the limit when $\Delta\energy\to 0$ of the
total number of eigenstates (normalized by area) with energies between
$\energy$ and $\energy+\Delta\energy$. The \emph{integrated}
density of states in turn is
$\IDoS(\energy)=\int_{x=-\infty}^\energy\DoS(x)\,dx$.
Let $\mathcal{F}$ and $\mathcal{B}$ denote the fundamental cell and
the first Brillouin zone respectively, and let
$E_{n}:\mathcal{B}\to\rset$ denote the $n$:th band in the 
bandstructure, that is $E_{n}(k)$ is the $n$:th smallest eigenvalue of
the algebraic eigenvalue problem~\eqref{eq:gen_evp} for $k\in\mathcal{B}$.
Then
\begin{align}
  \label{eq:IDoS}
  \IDoS(\energy) & = \frac{1}{|\mathcal{F}|} \sum_{n}
  \frac{1}{|\mathcal{B}|}\int_{\mathcal{B}}\chi_{\left\{\cdot<\energy\right\}}(E_{n}(k))\,dk,
\end{align}
where $\chi_{\left\{\cdot<\energy\right\}}$ is the indicator function
on the semi-infinite interval $(-\infty,\energy)$ and $|\cdot|$
denotes area.

The bands in~\eqref{eq:IDoS} are, 
in the case of the unperturbed graphene sheet on its fundamental
cell, 
$n\in\{1,2\}$ 
and for an $\supsz$-by-$\supsz$
super cell without vacancies $n\in\{1,2,\dots,2\supsz^2\}$. 
Similarly for the $\MoS_2$ model $n\in\{1,2,\dots,11\}$ and
$n\in\{1,2,\dots,11\supsz^2\}$, respectively.

For each sampled outcome~\eqref{eq:IDoS} is approximated from the
computed discretization 
of the bandstructure, $\{E_n(k_{lm})\}$ in two steps. 
First,  $E_n(k)$ is approximated by $\overline{E_{n}}(k)= E_n(k_{lm})$
where $k_{lm}$ is the discretization point closest to $k$.
%
Then, the indicator function~\eqref{eq:IDoS}
is approximated by a smoothed, Lipschitz continuous, step function
\begin{align}
  \label{eq:MLMC_smoothing}
  \chi_{\left\{\cdot<\energy\right\}}(E) & \approx
  g\left(\frac{E-\energy}{\smooth}\right), 
  \intertext{satisfying}
  g(x) & = 1, && \text{if $x \leq -1$},\nonumber \\
  g(x) & = 0, && \text{if $x \geq 1$},\nonumber \\
  \int_{-1}^1 
  x^q\left(\chi_{\left\{\cdot<0\right\}}(x)-g(x)\right)\,dx & = 0, 
  && \text{for $q=0,1$.}\nonumber 
\end{align}
This smoothing, using $\smooth\propto\tol$ where $\tol$ is the desired
accuracy, is needed when MLMC methods are used to compute distribution
functions of random variables; see~\cite{giles_distribution} for an
analysis of MLMC methods in this case. 
Similar smoothing strategies are also used in the computational
physics community.
Finally, $\IDoS$, is approximated in a uniform discretization
$\energy_0<\energy_1<\dots<\energy_M$ of an interval 
containing the range of computed energies. 

The expected value of the integrated density of states is approximated by
Monte Carlo or MLMC sample averages. From the expected value of the
integrated density of states the density of states may be estimated by
numerical differentiation.


\subsection{Numerical Results}
\label{sec:numres}

The following numerical results are intended to show whether an MLMC
approach can accelerate computations of the quantity of interest in
the test problems; in particular it is important to see that the
control variates suggested in Section~\ref{sec:CV} improves on the
rate of convergence of the variance of the samples, so that $D>S$ in
the models~\eqref{eq:compl_model} and~\eqref{eq:compl_model2}.

\paragraph{A Tight-binding model of graphene}

\begin{figure} \centering
  \includegraphics[width=0.45\textwidth]{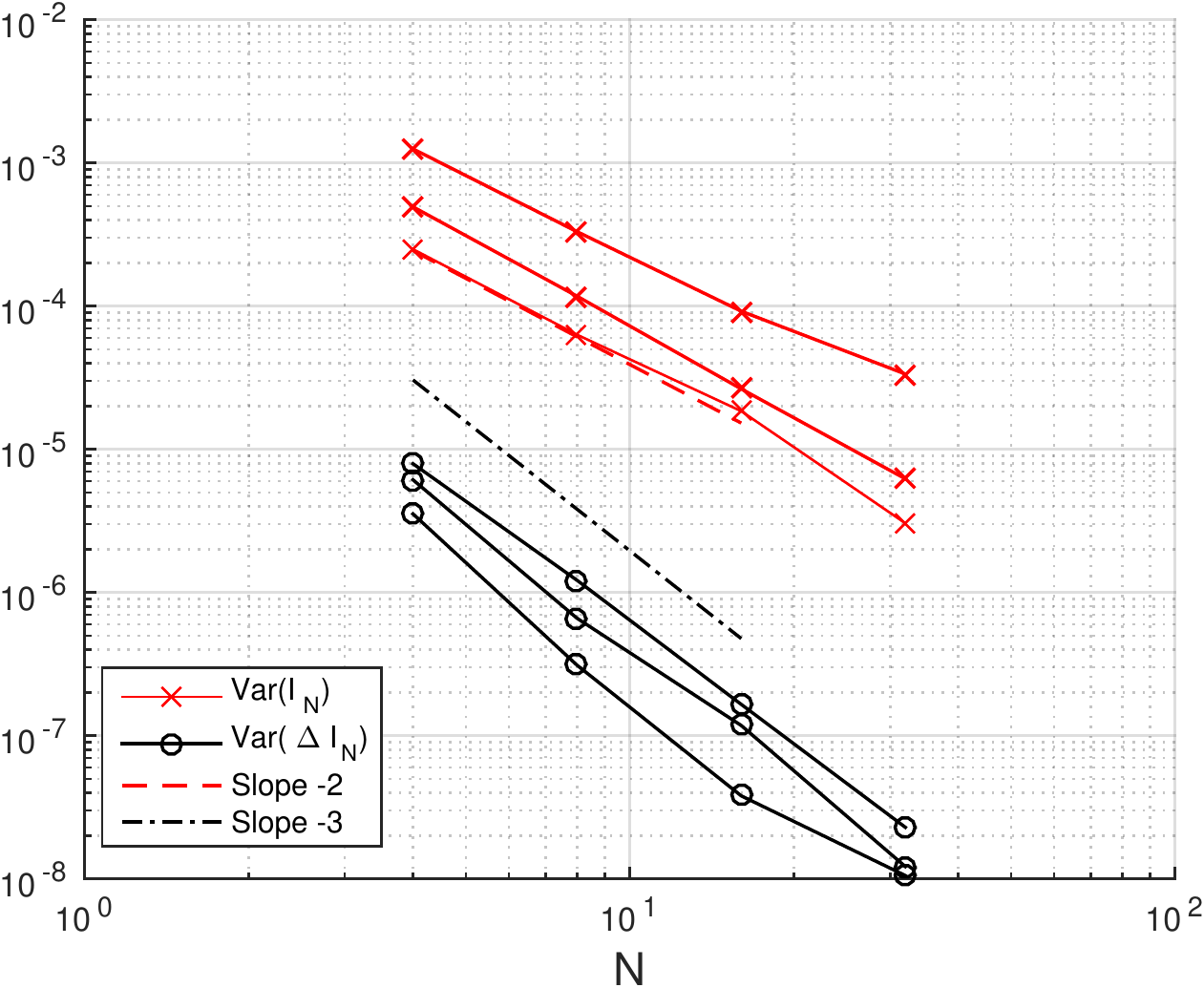}
  \includegraphics[width=0.45\textwidth]{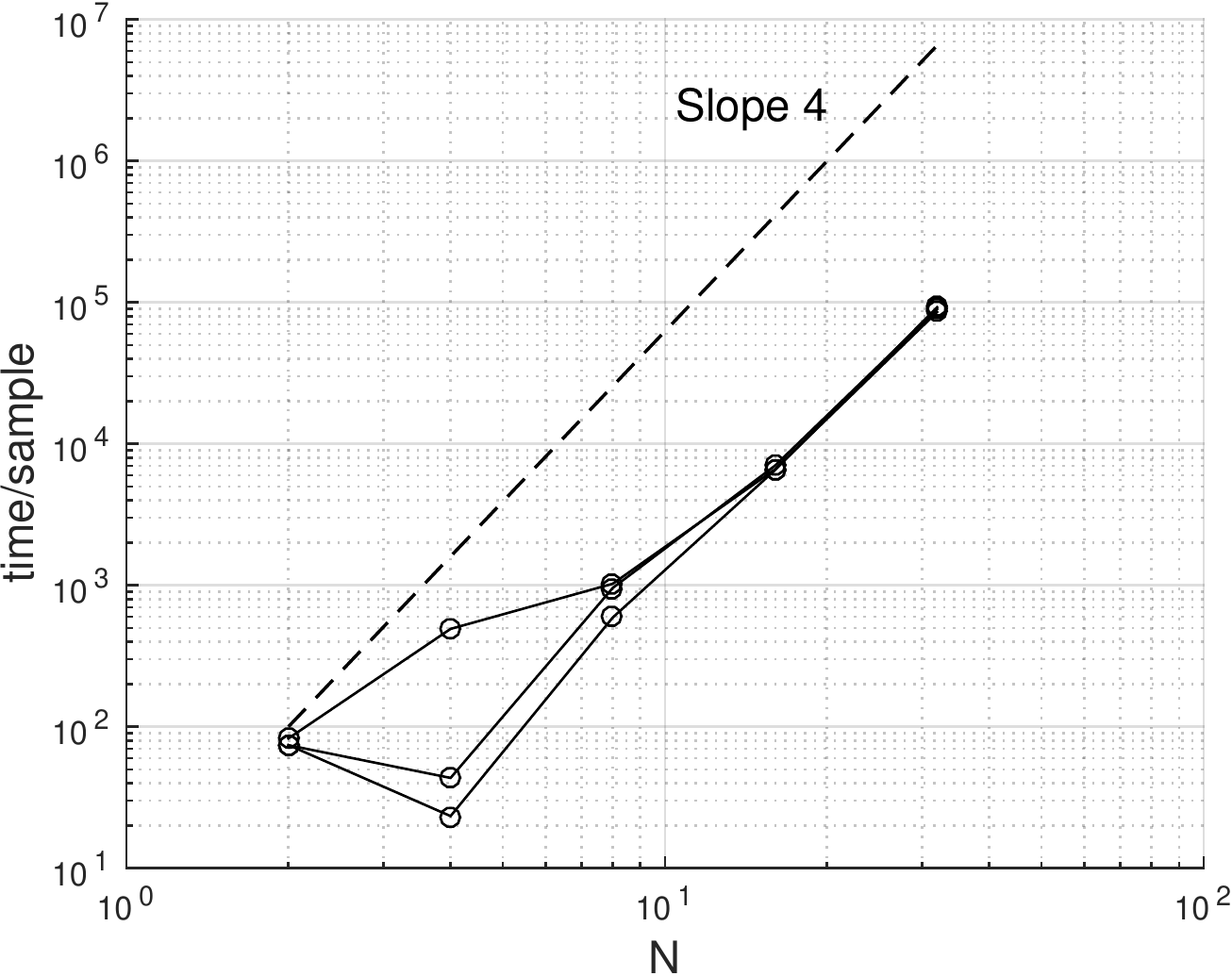}
  \caption{$\mathrm{MoS_2}$: (Left) The sample variance of the integrated
    density of states per unit area, $\IDoS_\supsz(\energy)$, using a
    super cell of size $\supsz\times\supsz$ and the sample variance
    of the difference, 
    $\Delta\IDoS_\supsz(\energy)=\IDoS_\supsz(\energy)-\IDoS_{\supsz/2}(\energy)$,
    for the three vacancy probabilities in
    Table~\ref{tab:MoS2_MLMC}. Shown here is the arithmetic mean of
    the quantities over the discretization points in the interval 
    $-6 \mathrm{eV}<\energy<4 \mathrm{eV}$ and the sample variance was
    computed using the samples in the MLMC estimators. In particular the
    sample variance on the largest super cell is based on only five
    samples. The experimentally observed convergence rates are
    approximately $S=2$ and $D=3$.\\
    (Right) Wall time per sample in the simulations where each sample
    was run on a single core of a multi core processor. An
    eigenvalue problem for a full matrix of side $\propto\supsz^2$
    were solved for every discretization point of the Brillouin
    zone, giving the cost per eigenvalue solve
    $\propto\supsz^6$ for large enough $\supsz$. Since the number of
    such discretization points were chosen to decrease as
    $\supsz^{-2}$, the observed time per sample is approximately
    $\propto\supsz^4$; see Section~\ref{sec:bandstructure}.} 
  \label{fig:Conv_MoS2}
\end{figure}

An empirical investigation of how the quantities used in the
complexity analysis of Section~\ref{sec:MC_MLMC} behave for the
tight-binding model of graphene using modest super cell sizes, up to a
32-by-32 extension of the fundamental cell of the honeycomb lattice,
containing 2048 atom locations. 
The results show that in this example the sample variance of the
quantity of interest, $Q_\ell$, measured in discrete norms, decays approximately
as $\supsz^{-2}$, and the sample variance of $Q_\ell-Q_\ell^{CV}$
decays faster, approximately as $\supsz^{-3}$. 
The computational cost per sample is nearly independent of
$\supsz$ for the first few sample points, where the generalized
eigenvalue problems only involve a few unknowns, and starts to grow
only around $\supsz=8$. Between $\supsz=16$ and $\supsz=32$ the rate
of growth is approximately 4; see also Figure~\ref{fig:Conv_MoS2} for
the $\MoS_2$ case.
In the notation of Section~\ref{sec:MC_MLMC}, the empirical estimates
of the parameters are
\begin{align}
  \label{eq:param_est}
  W & \approx 1.5, & S & = 2, & D & = 3, & C & = 4.
\end{align}
Since $D>S$ the asymptotic complexity of an MLMC algorithm should be
better than that of a standard Monte Carlo method. 
We expect 
an improvement on the computational
work using MLMC as soon as $\supsz>=32$ here. The smallest control variate
worth including in the MLMC estimator~\eqref{eq:MLMC} is $\supsz=16$ since
samples on smaller super cell sizes are nearly as expensive.

\begin{figure}
  \centering
  \includegraphics[width=0.45\textwidth]{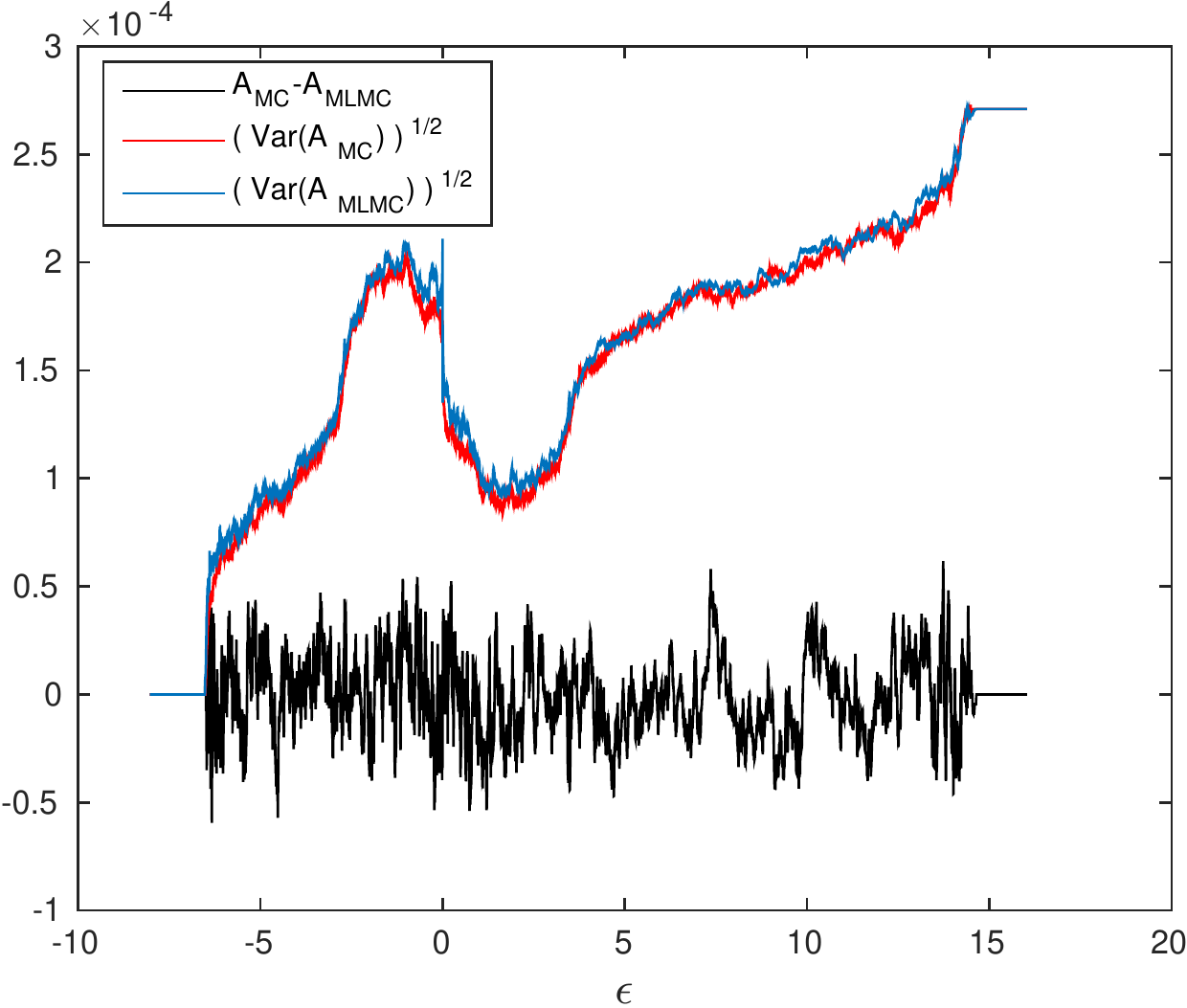}
  \includegraphics[width=0.45\textwidth]{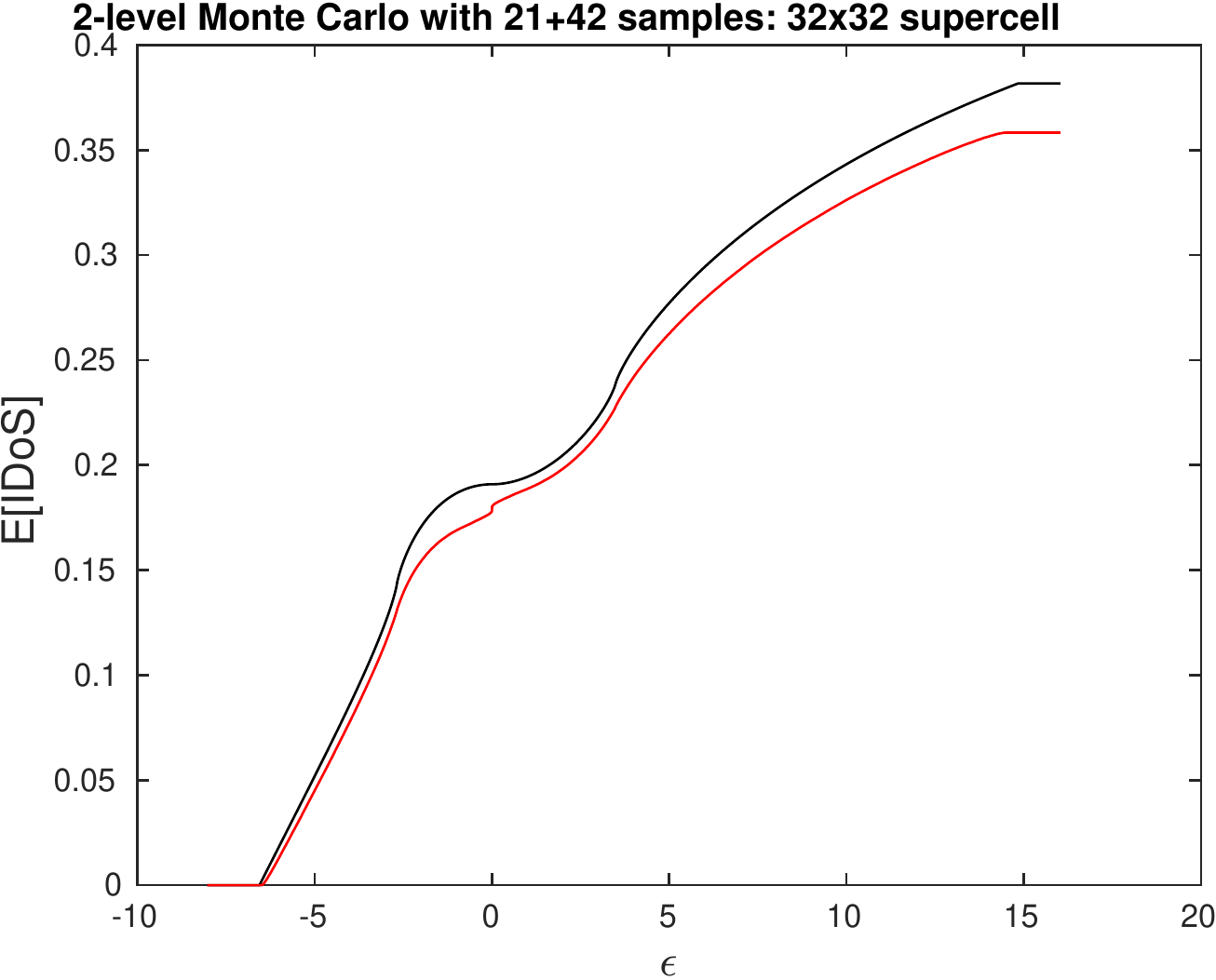}
  \includegraphics[width=0.45\textwidth]{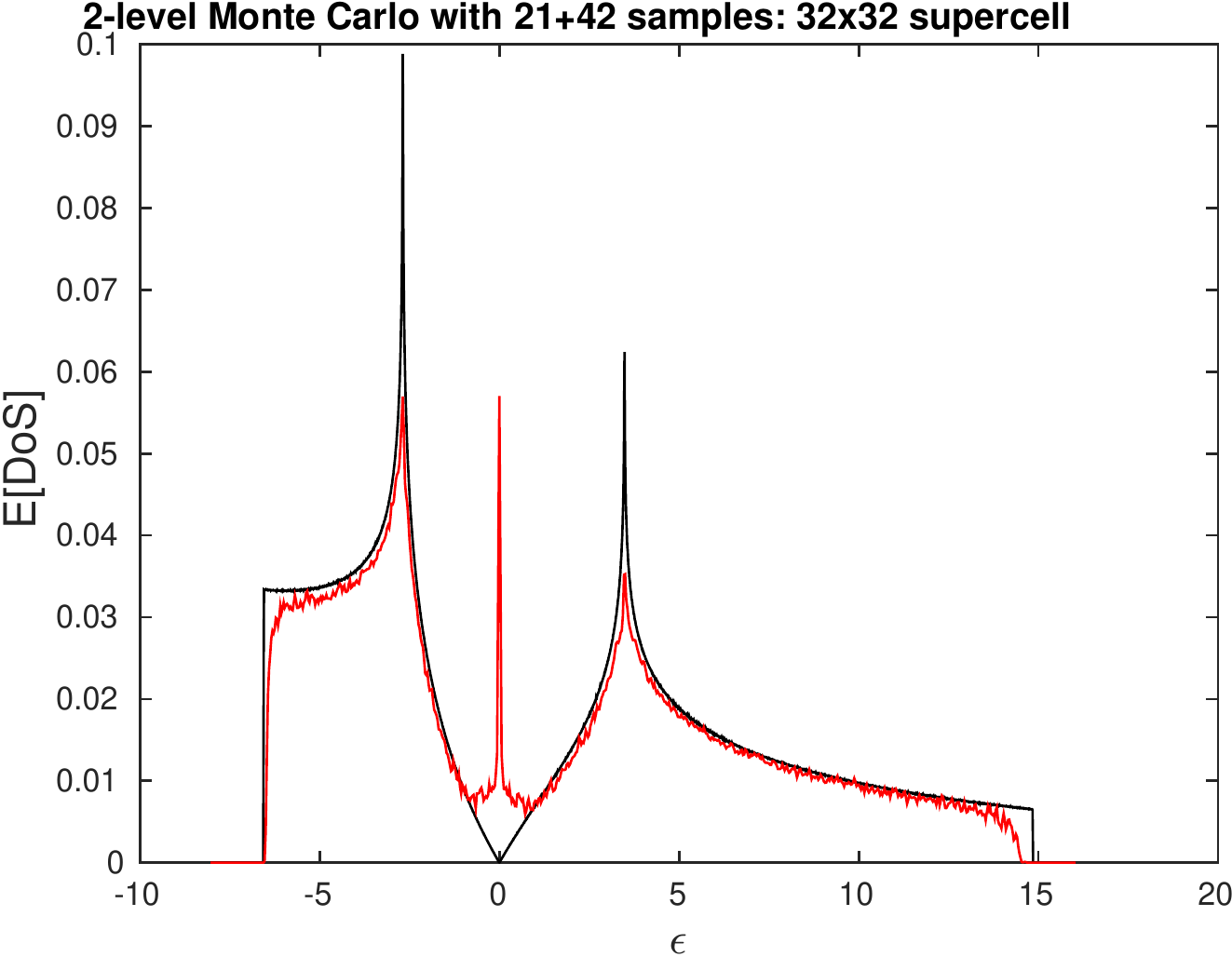}
  \includegraphics[width=0.45\textwidth]{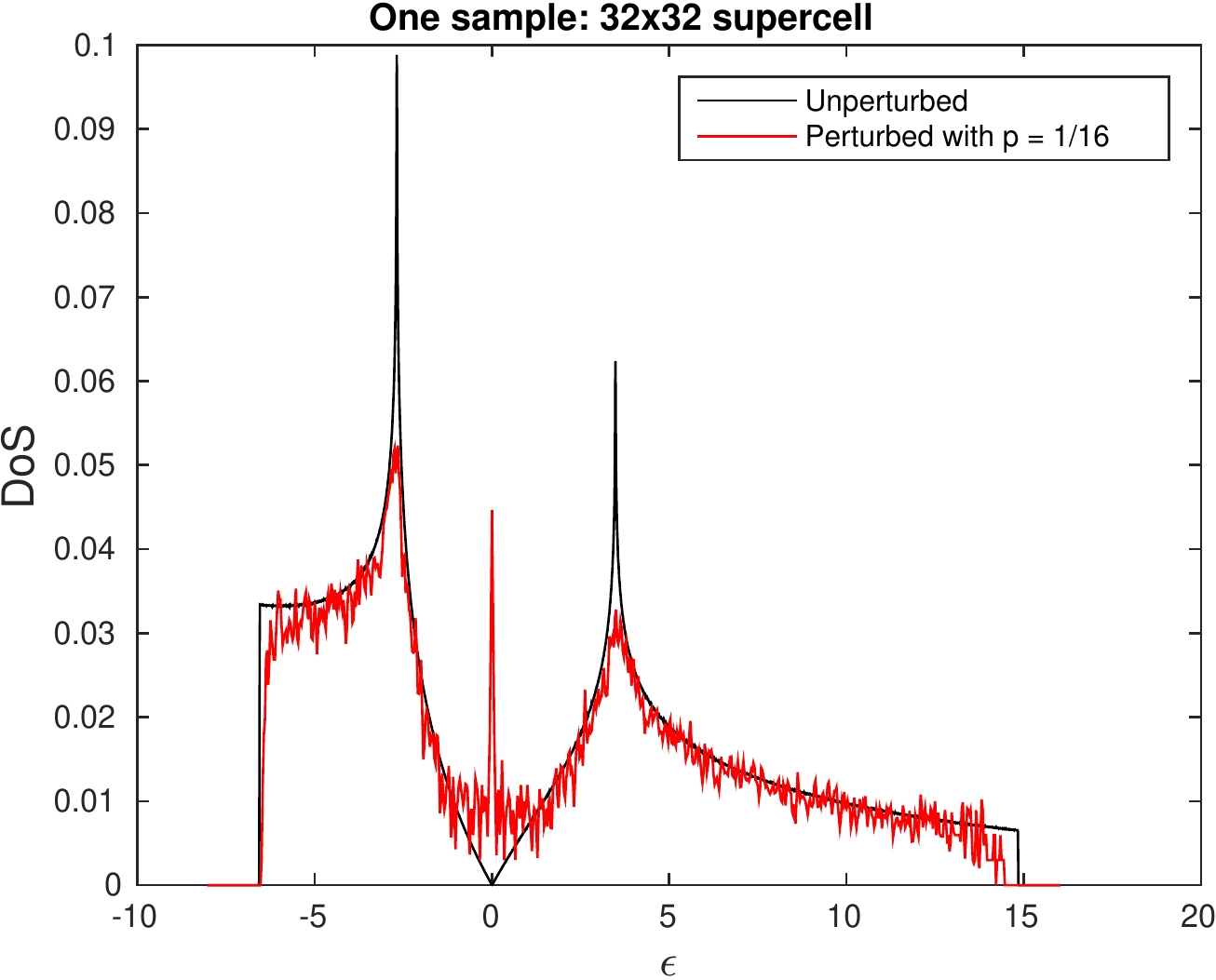}
  \caption{A bi-level Monte Carlo approximation of the integrated
    density of states on a 32-by-32 super cell with
    probability of any atom location being
    removed from the tight-binding model $\Pvac=0.0625$, denoted $Q_\ell$ below.\\
    (Top left) Black curve shows the difference between a 42 sample
    standard Monte Carlo estimate of $Q_\ell$ and a bi-level Monte Carlo
    estimator using 21 samples of $Q_\ell$ and 42 of its control variate
    $Q_\ell^{CV}$, obtained at $61\%$ of the cost of the single
    level. The standard deviations of the two estimators are of the same order.\\
    (Top right) The bi-level Monte Carlo estimate of $\E{Q_\ell}$
    together with the unperturbed.\\
    (Bottom left) Approximation of the density of states obtained by
    numerical differentiation of the
    bi-level Monte Carlo estimate above.\\
    (Bottom right) Approximation of the density of states based on
    only one sample and the same resolution in the energy. 
  }
  \label{fig:BiLMC}
\end{figure}

Following the observation above, a 2-level Monte Carlo estimator based
on super cell sizes $\supsz=32$ and $\supsz=16$ for the control variate is shown
in Figure~\ref{fig:BiLMC}. Here the 2-level estimator used 21 samples
on the larger super cell size, $\supsz=32$, and 42 samples on the smaller
size, $\supsz=16$. For comparison an additional 21 independent samples on
$\supsz=32$ were generated and a single level estimator based on 42 samples
computed. The variance of the two estimators are nearly of the same
magnitude as desired, while the cost of the 2-level estimator was
$61\%$ of that of the standard Monte Carlo estimator. 
It can be seen most clearly from the density of states, computed by
numerical differentiation, that it is crucial to control the
statistical error even on a super cell of this size. The two plots of
the density of states computed either from the 2-level Monte Carlo
estimator or from a single outcome of random impurities use the same
resolution in the energy; in the latter case noise hides all detail.

Note that the work ratio between MLMC and standard Monte Carlo will
not remain constant at around $61\%$ as we aim for more accurate
solutions, provided that the empirical complexity and convergence
estimates extrapolate to larger $\supsz$ with the present rates. The
next example will illustrate this.

\paragraph{The tight-binding model of $\MoS_2$}

Here, using the parameters in
Table~\ref{tab:MoS2_MLMC}, we again observe the values
in~\eqref{eq:param_est} for the parameters in the convergence and work
models. By the estimates of Section~\ref{sec:MC_MLMC}, we expect the
computation time of standard Monte Carlo to grow as
$\tol^{-(2+\frac{C-S}{W})}\approx\tol^{-10/3}$ while that of MLMC to
grow as $\tol^{-(2+\frac{C-D}{W})}\approx\tol^{-8/3}$ as $\tol\to 0$.
For the fixed accuracy of the numerical results here, 
we estimate that using a standard Monte Carlo estimator the work
required to obtain a variance in $\IDoS(\energy)$ comparable to that
observed in the MLMC estimator would be one order of magnitude larger;
see Table~\ref{tab:MoS2_conv}.

In the numerical tests of the $\MoS_2$ we made use of the fact
that in the total number of possible permutations is finite for each
finite super cell size. 
For a sufficiently small super cell, the number of possible
combinations is small enough to compute the quantity of interest on
all of them, taking symmetries into account, and then by combinatorial
means obtain the probabilities of all possible outcomes for a complete
description of the statistics. 
This was done for the smallest $2\times2$ super cell for all values of
$\Pvac$. 
For $\Pvac=0.025$ and 0.05 we also took advantage of the finite space of
possible outcomes of the perturbations by identifying identical
samples of the random perturbations beforehand and avoiding repeating
them. This leads to substantial computational savings on the still
rather small $4\times4$ and $8\times8$ super cells.

In these numerical tests we started with rough estimates of
the parameters in the models~\eqref{eq:compl_model}
and~\eqref{eq:compl_model2} to determine a sequence of samples
using~\eqref{eq:opt_samples}. The resulting values of $\nrsam_\ell$
are rough approximations of the optimal choices. An alternative
approach is to use an algorithm to estimate the parameters during the
computation and adaptively choose the number of samples;
see~\cite{Collier_CMLMC}.

\begin{table}
  \centering
  \begin{tabular}{|c|c|c|c|c||c||c|}
    \hline
    \multicolumn{7}{c}{$\Pvac=0.025$}\\
    \hline
    Level & $\supsz$ & $\nrsam$ & $\supsz\nrBZp$ & $\smooth$ &
    $\Delta\energy$  & time (h)\\ 
    \hline
    1 & 2 & Exhaustive & 128 & 0.01 & 3.9\e{-3} & 0.33 \\
    \hline
    2 & 4  & 2072 & 128 & 0.01 & 3.9\e{-3} & 13.4 \\
    3 & 8  & 564  & 128 & 0.01 & 3.9\e{-3} & 92.7 \\
    4 & 16 & 76   & 128 & 0.01 & 3.9\e{-3} & 137  \\
    5 & 32 & 5    & 128 & 0.01 & 3.9\e{-3} & 128  \\    
    \hline
     & & & & & total time & 372\\
    \hline
    \multicolumn{7}{c}{ }\\
    \hline
    \multicolumn{7}{c}{$\Pvac=0.05$}\\
    \hline
    Level & $\supsz$ & $\nrsam$ & $\supsz\nrBZp$ & $\smooth$ &
    $\Delta\energy$  & time (h)\\ 
    \hline
    1 & 2 & Exhaustive & 128 & 0.01 & 7.8\e{-3} & 0.33 \\
    \hline
    2 & 4  & 2450 & 128 & 0.01 & 7.8\e{-3} & 30 \\
    3 & 8  & 474  & 128 & 0.01 & 7.8\e{-3} & 126 \\
    4 & 16 & 77   & 128 & 0.01 & 7.8\e{-3} & 150 \\
    5 & 32 & 5    & 128 & 0.01 & 7.8\e{-3} & 124 \\    
    \hline
     & & & & & total time & 430\\
    \hline
    \multicolumn{7}{c}{ }\\
    \hline
    \multicolumn{7}{c}{$\Pvac=0.1$}\\
    \hline
    Level & $\supsz$ & $\nrsam$ & $\supsz\nrBZp$ & $\smooth$ &
    $\Delta\energy$  & time (h)\\ 
    \hline
    1 & 2 & Exhaustive & 128 & 0.01 & 15.6\e{-3} & 0.37 \\
    \hline
    2 & 4  & 2072 & 128 & 0.01 & 15.6\e{-3} & 283 \\
    3 & 8  & 564  & 128 & 0.01 & 15.6\e{-3} & 161 \\
    4 & 16 & 76   & 128 & 0.01 & 15.6\e{-3} & 139  \\
    5 & 32 & 5    & 128 & 0.01 & 15.6\e{-3} & 129  \\    
    \hline
     & & & & & total time & 702\\
    \hline
  \end{tabular}
  \caption{Parameters in the MLMC estimator in
    Figure~\ref{fig:IDoS_MoS2} and the computational times spent on
    each level of the MLMC hierarchy as well as the total time. Here, 
    $\supsz$ is the super cell size, 
    $\nrsam$ is the number of samples, 
    $\nrBZp$ controls the discretization of the Brillouin zone as 
      in Section~\ref{sec:bandstructure}, 
    $\smooth$ is the smoothing parameter
      in~\eqref{eq:MLMC_smoothing}, and 
    $\Delta\energy$ is the step size in the numerical
      differentiation in the post processing step used to get the
      density of states in Figure~\ref{fig:DoS_MoS2}.\\
    The computational times are wall times for one core on multi-core
    processors, where one sample was running on each core. 
    The computations with $\Pvac=0.025$ and  $\Pvac=0.05$ did not repeat
    computations on identical outcomes of the random perturbation
    leading to significant computational savings on levels 1 and 2
    where the probability of repeated outcomes is high.
  }
  \label{tab:MoS2_MLMC}
\end{table}

\begin{table}
  \centering
  \begin{tabular}{|c||c|c|c|c|c|c|c|c|}
    \hline
    $\Pvac$ & $W$ & $S$ & $D$ & $C$ & 
    $AC_{\mathtt{FS}}$  & $AC_{\mathtt{SLMC}}$  & $AC_{\mathtt{MLMC}}$ & $R$ \\
    \hline
    0.025 & 3/2 & 2 & 3 & 4 & 4 & 3+1/3 & 2+2/3 & 0.06 \\
    0.05  & 3/2 & 2 & 3 & 4 & 4 & 3+1/3 & 2+2/3 & 0.07 \\
    0.1   & 3/2 & 2 & 3 & 4 & 4 & 3+1/3 & 2+2/3 & 0.06 \\
    \hline
  \end{tabular}
  \caption{The parameters, $W$, $S$, $D$, and $C$, in the
    models~\eqref{eq:compl_model} and~\eqref{eq:compl_model2}
    estimated from the numerical experiment on $\MoS_2$;
    compare Figure~\ref{fig:Conv_MoS2}. Included are also the
    corresponding estimated asymptotic complexities
    $\mathtt{Work}\propto\tol^{-AC}$ in the work
    estimates~\eqref{eq:Work_FS},~\eqref{eq:Work_SLMC} and~\eqref{eq:Work_MLMC}.
    Finally, $R$ denotes the ratio between the observed computational
    time of the MLMC method and the \emph{estimated} time for a
    standard Monte Carlo method to obtain approximately the same
    variance; see Figure~\ref{fig:Var_IDoS_MoS2}.}
  \label{tab:MoS2_conv}
\end{table}

\begin{figure} \centering
  \includegraphics[width=0.45\textwidth]{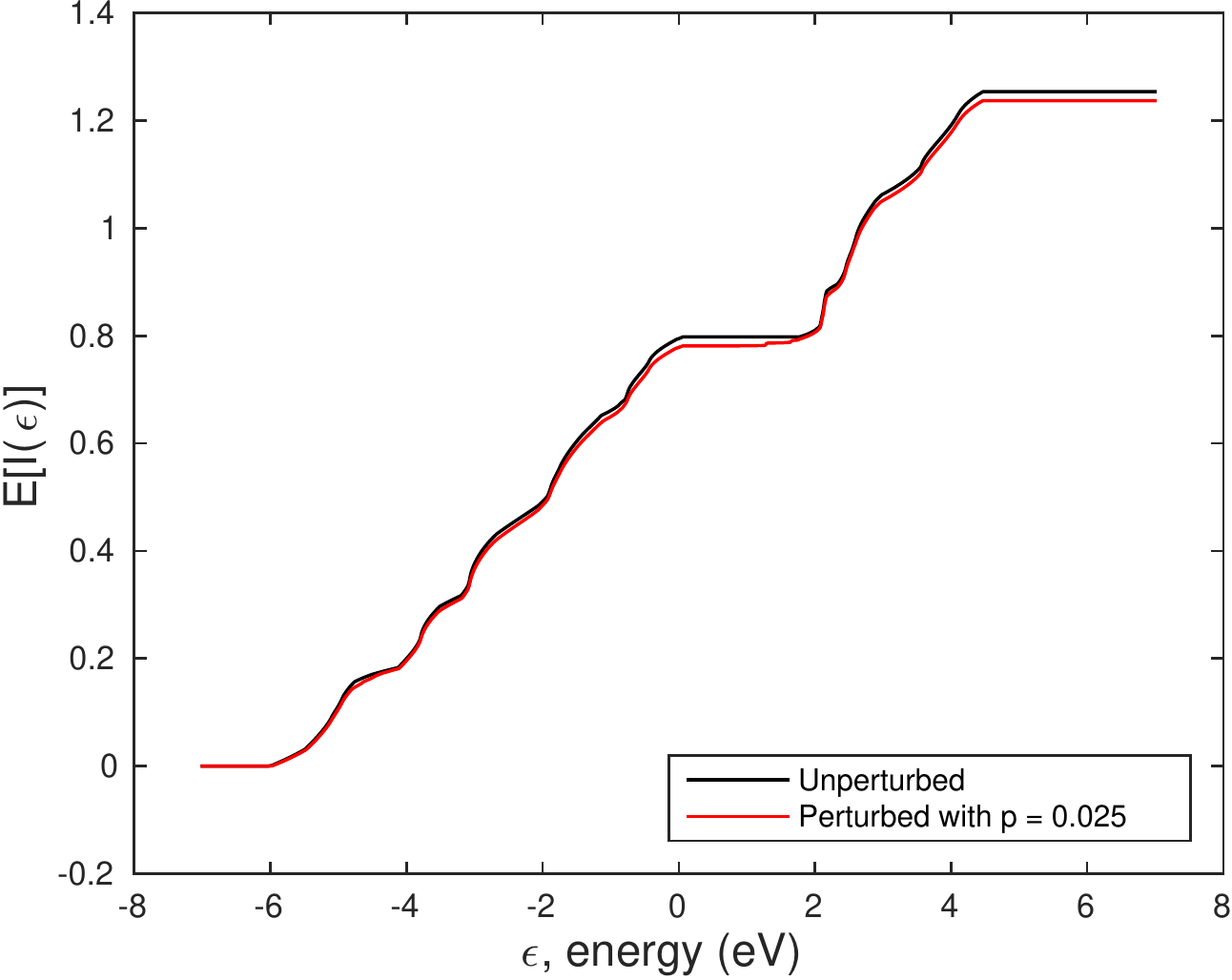}
  \includegraphics[width=0.45\textwidth]{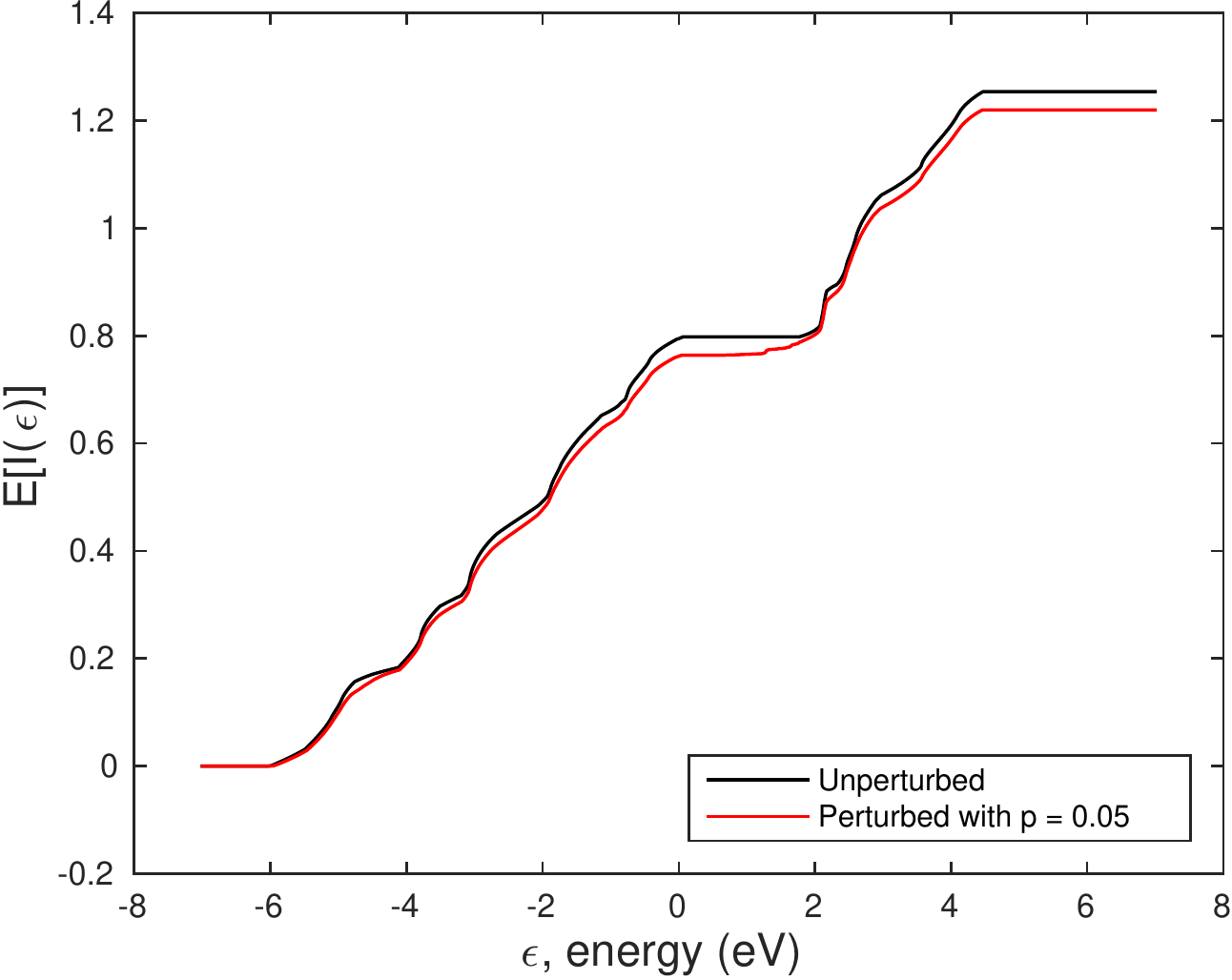}
  \includegraphics[width=0.45\textwidth]{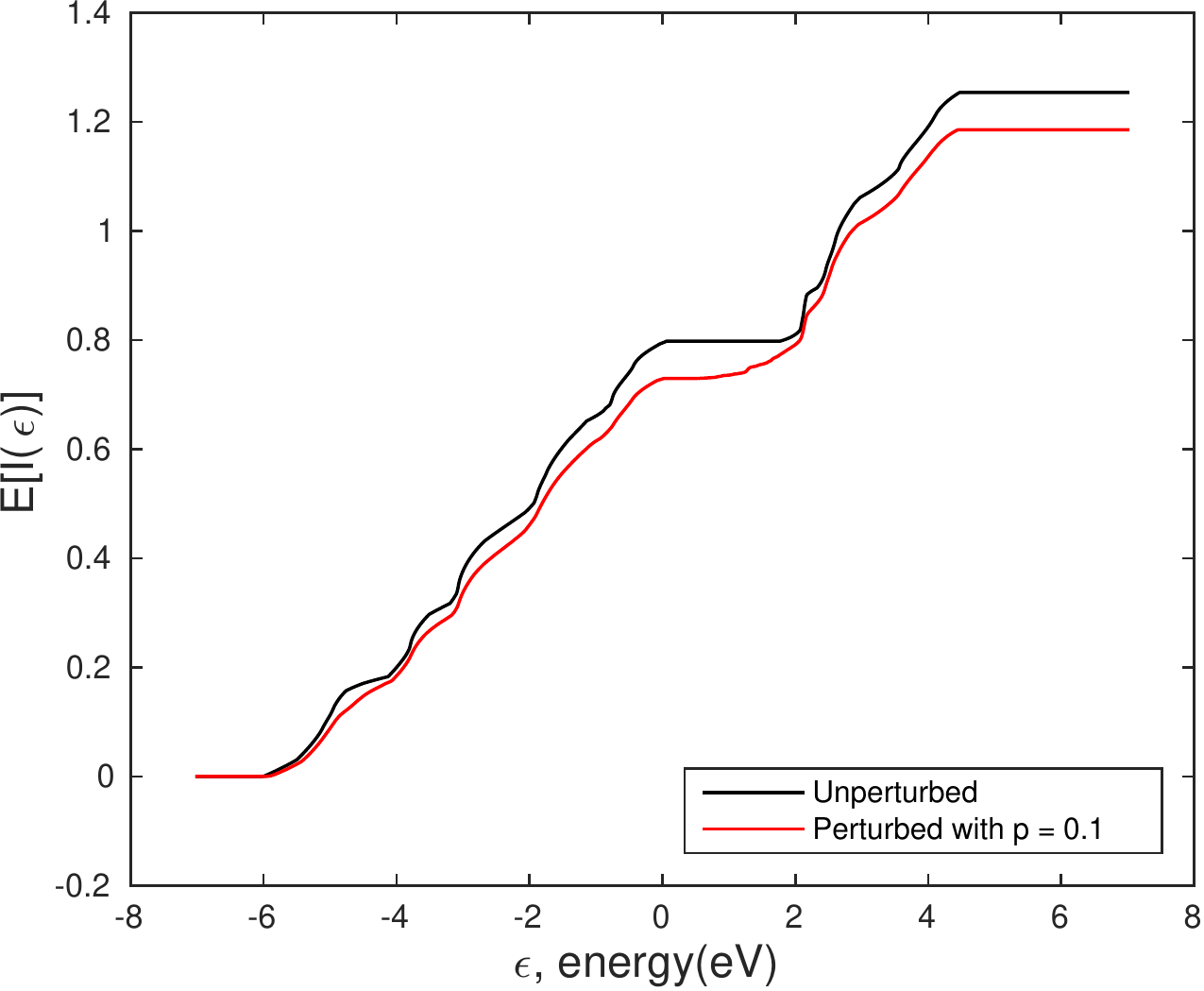}
  \caption{$\mathrm{MoS_2}$: 
    MLMC approximations of the expected integrated
    density of states, $\IDoS(\epsilon)$, on a 32-by-32 super cell
    with the probability, $\Pvac$, of any $\mathrm{S}$ atom pair 
    being removed from the tight-binding model taking the values
    $\Pvac=0.025,~0.05~,0.1$ respectively. 
    The integrated density of states for unperturbed material
    is shown for comparison. 
    The MLMC estimators were computed using the parameters in
    Table~\ref{tab:MoS2_MLMC}.}
  \label{fig:IDoS_MoS2}
\end{figure}

\begin{figure} \centering
  \includegraphics[width=0.45\textwidth]{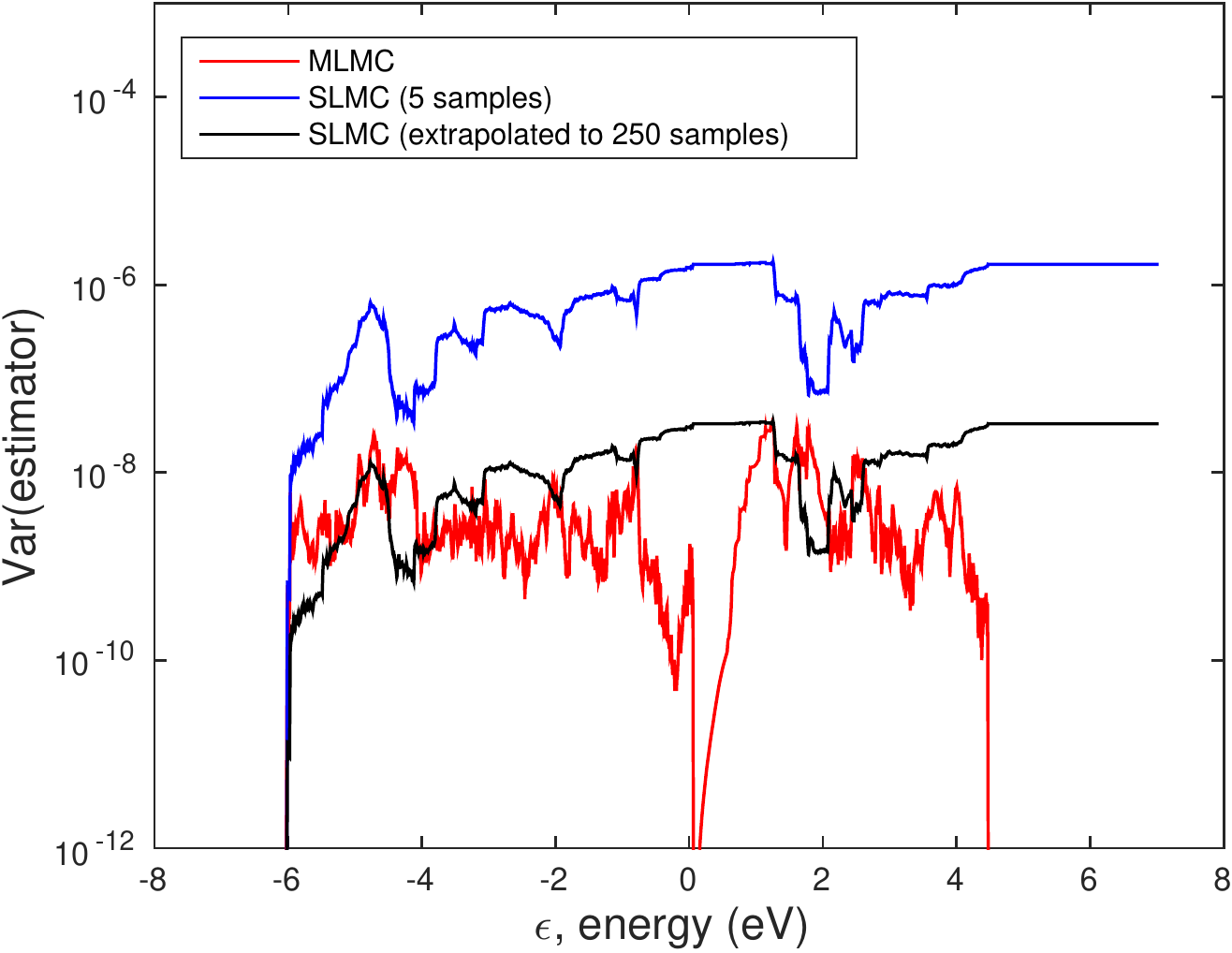}
  \includegraphics[width=0.45\textwidth]{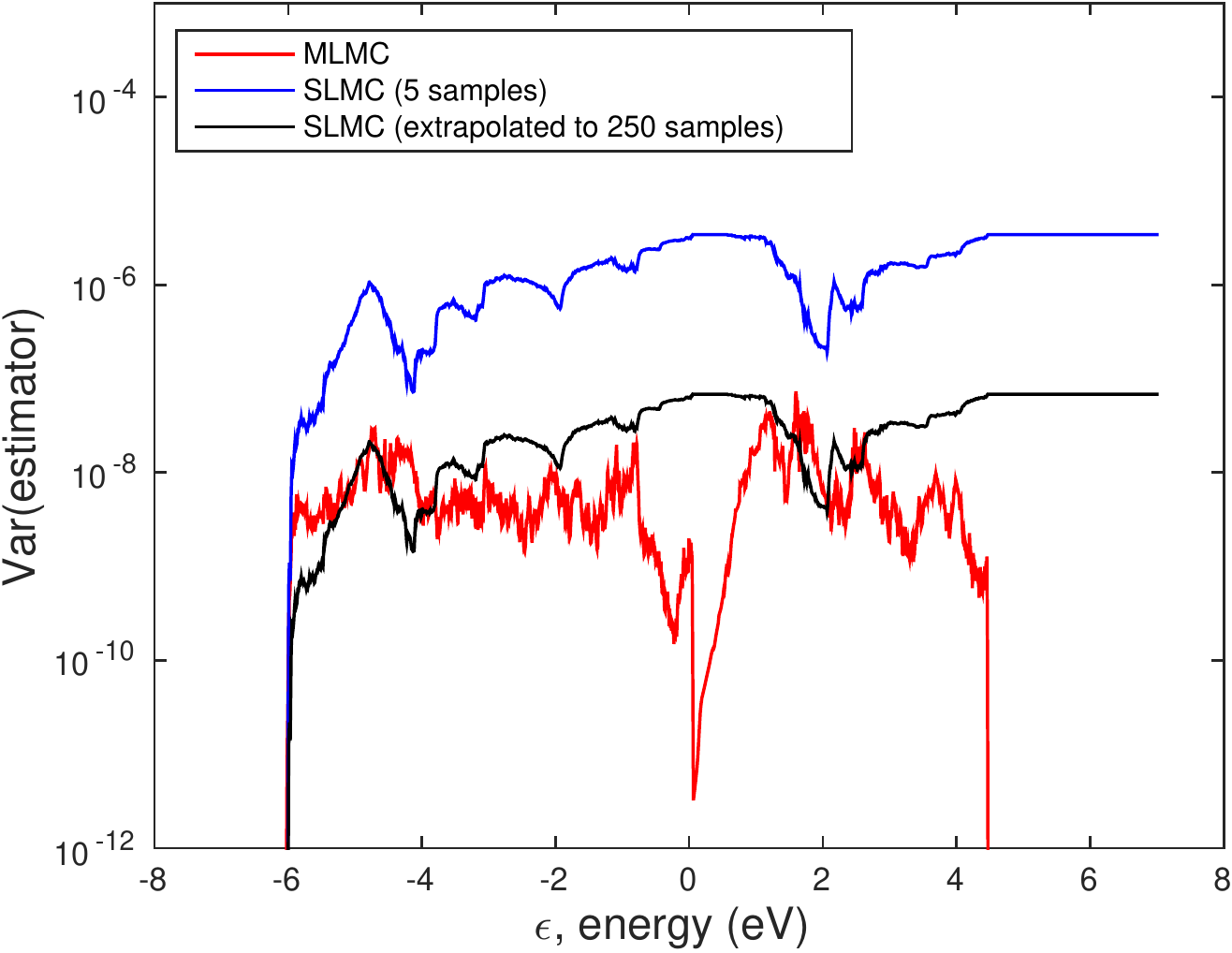}
  \includegraphics[width=0.45\textwidth]{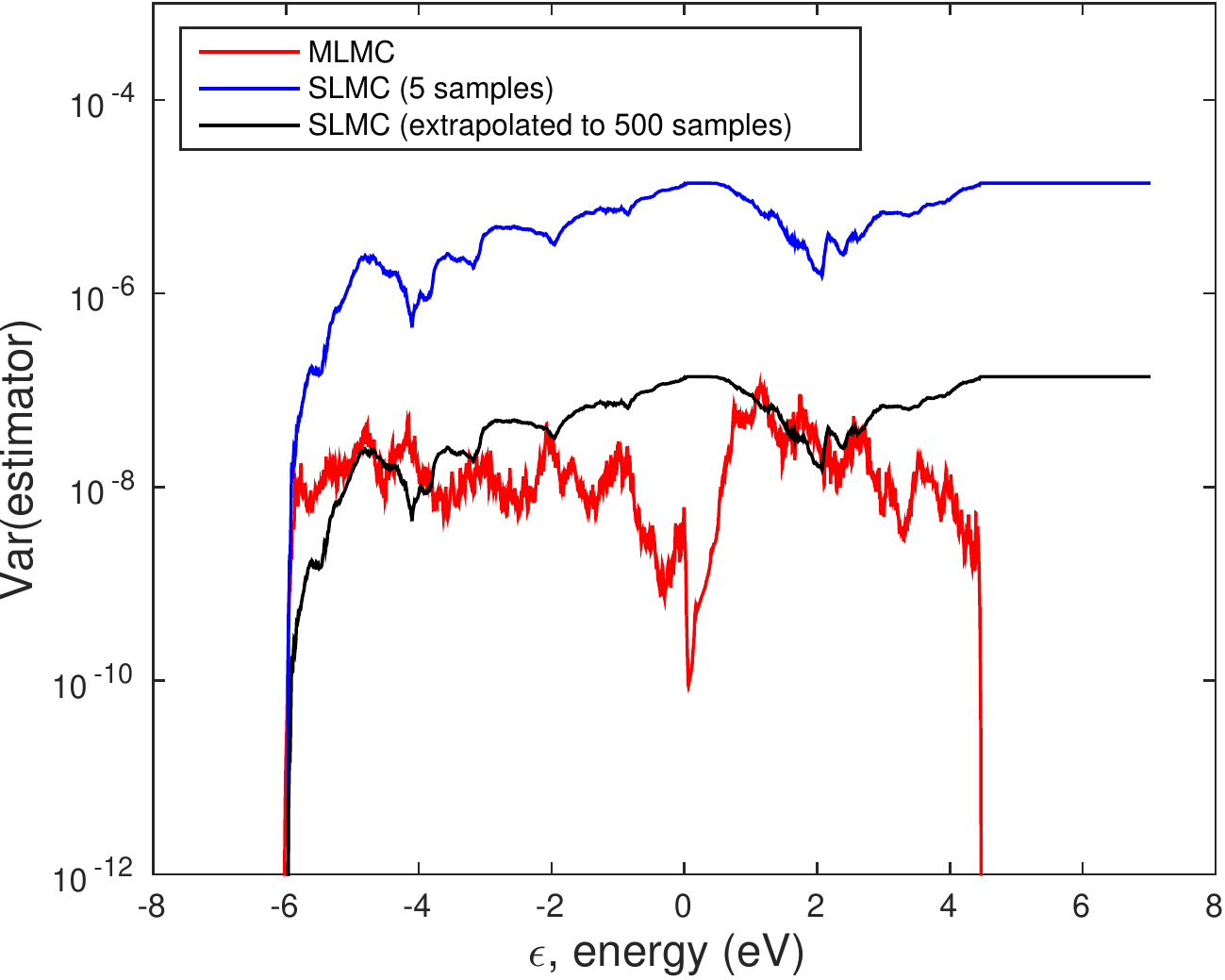}
  \caption{$\mathrm{MoS_2}$: 
    Estimates of the pointwise variance of the MLMC estimators of
    Figure~\ref{fig:IDoS_MoS2} are compared with the corresponding
    variance estimates when only the five samples on the 32-by-32 super cell
    were used in a single level Monte Carlo (SLMC) estimator. 
    Also included are rescaled versions of the SLMC variances chosen
    so that they are comparable to those of the MLMC
    estimators in the interesting range $1 \mathrm{eV}<\energy<2
    \mathrm{eV}$, which contains the upper part of the band gap of the
    unperturbed material.
    This gives rough estimates of how many samples the SLMC estimators
    would need to match the error of the MLMC estimators; see
    Table~\ref{tab:MoS2_MLMC}.} 
  \label{fig:Var_IDoS_MoS2}
\end{figure}

\begin{figure} \centering
  \includegraphics[width=0.45\textwidth]{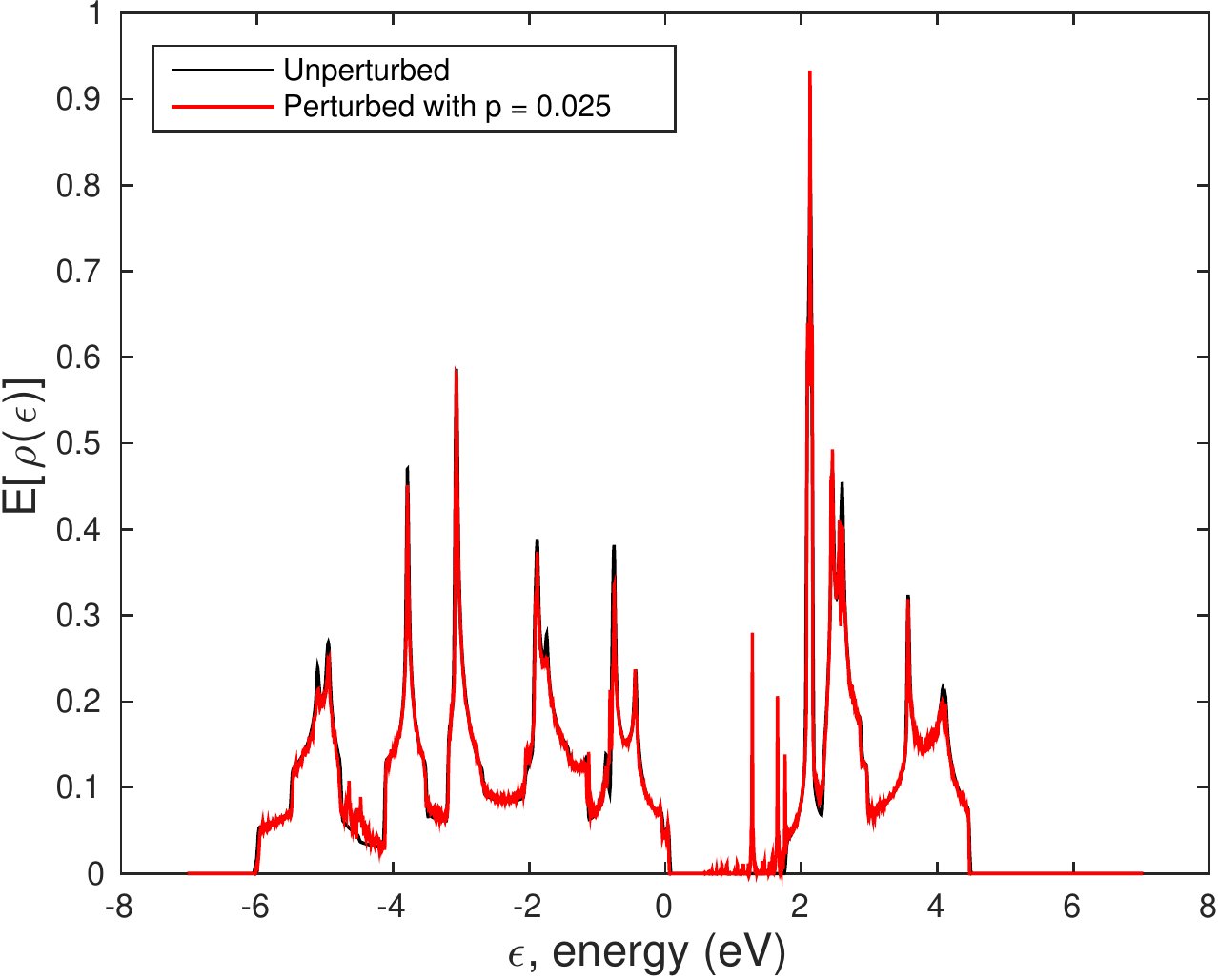}
  \includegraphics[width=0.45\textwidth]{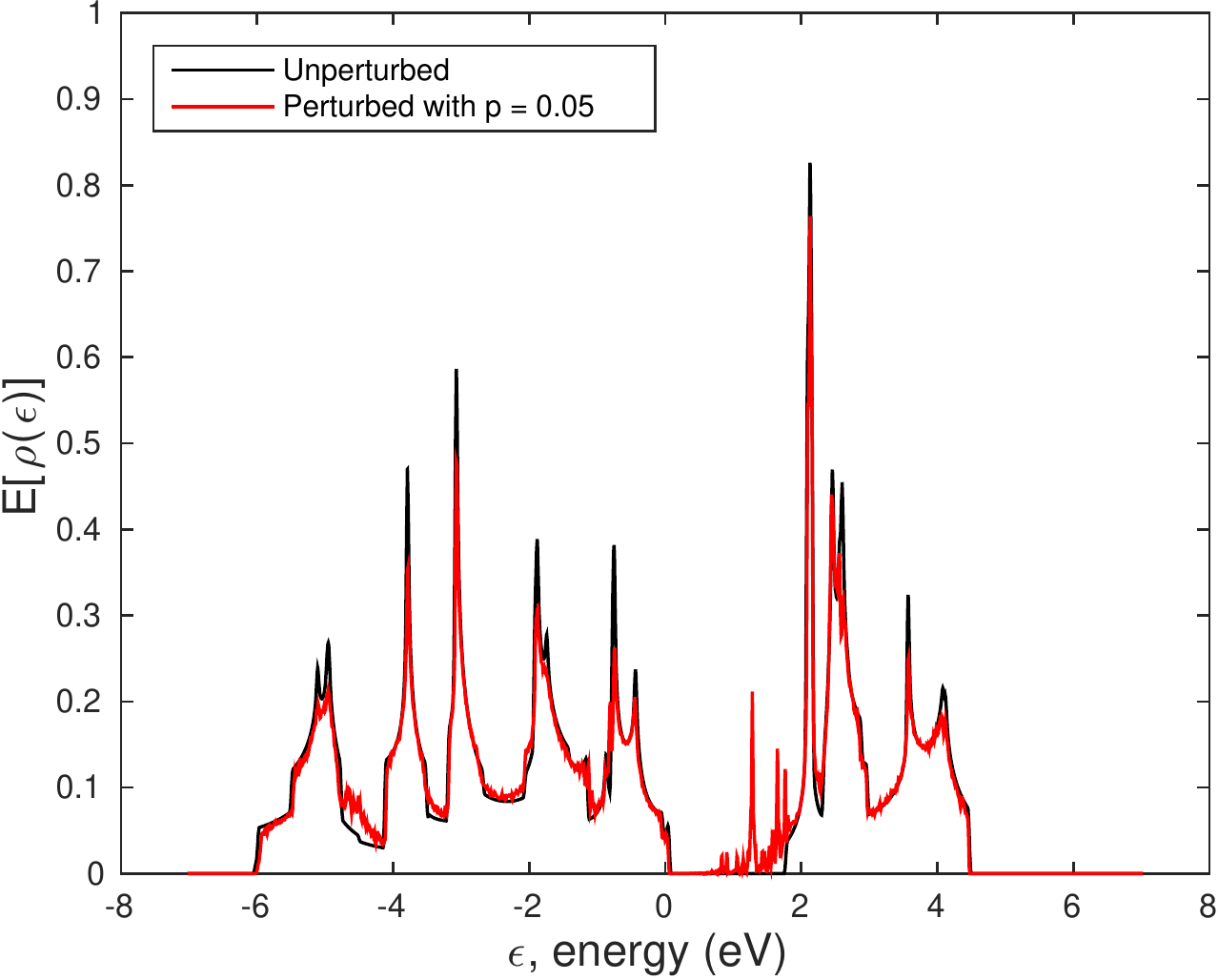}
  \includegraphics[width=0.45\textwidth]{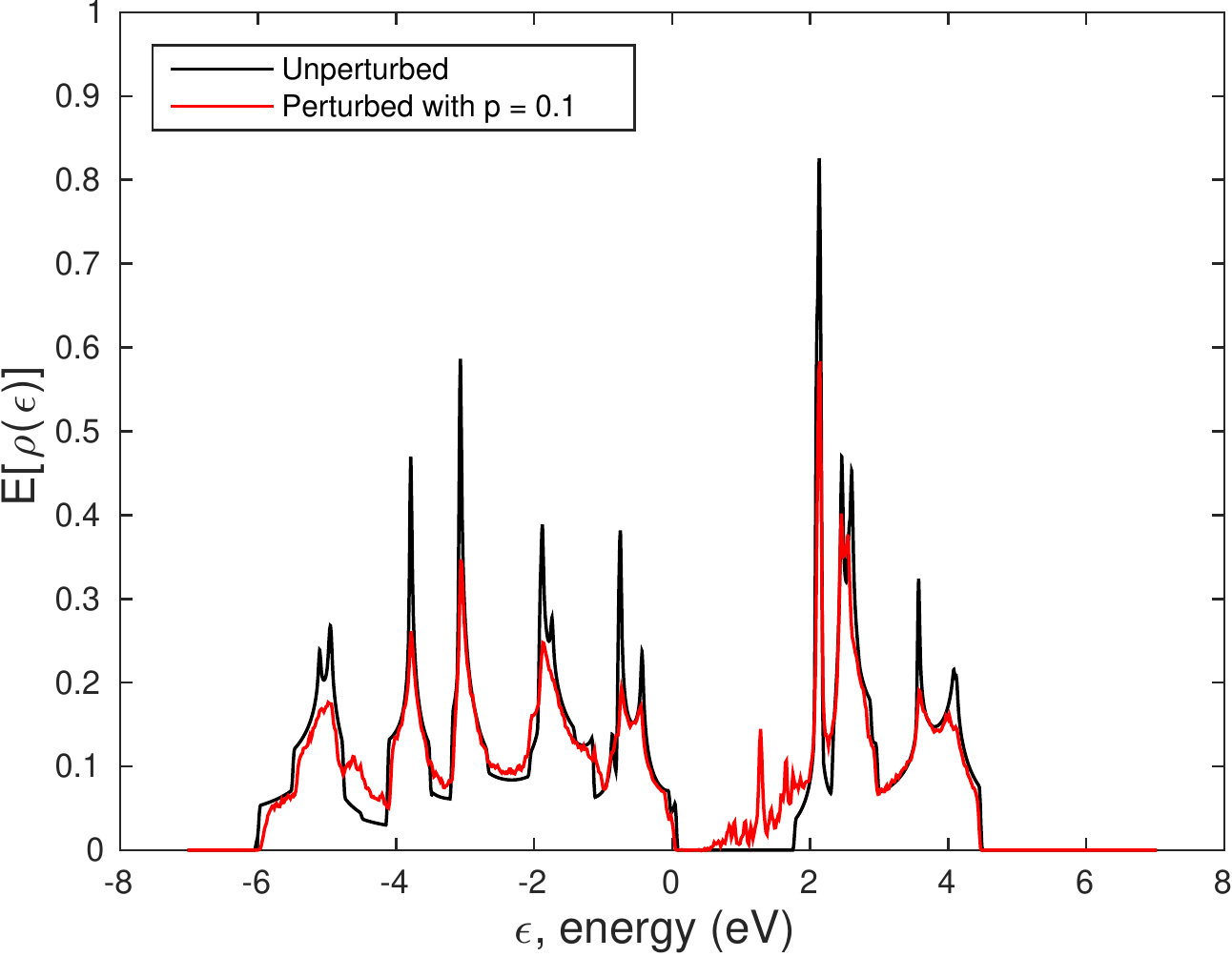}
  \includegraphics[width=0.45\textwidth]{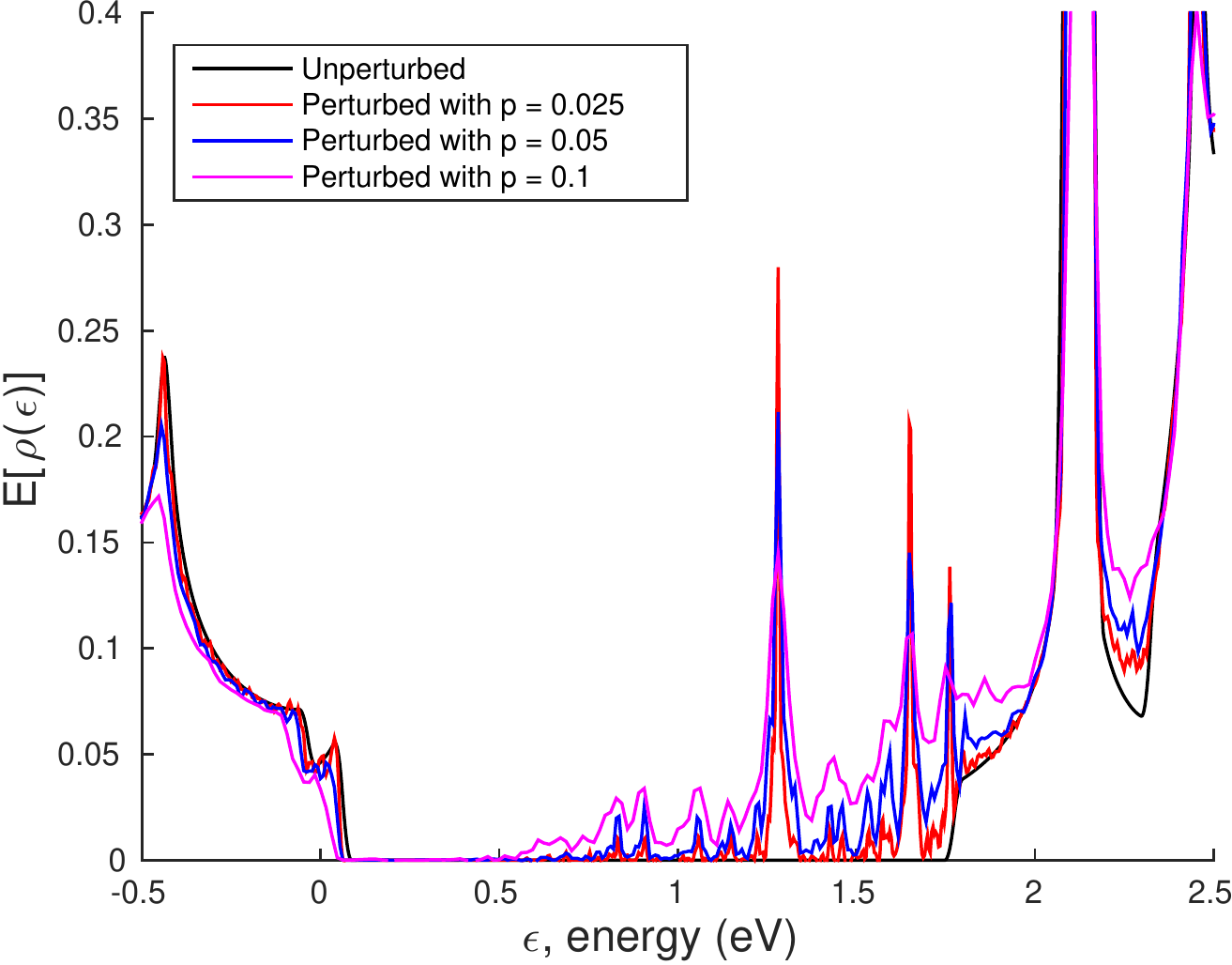}
  \caption{$\mathrm{MoS_2}$: The density of states per unit area,
    $\DoS(\energy)$, computed by numerical differentiation of the MLMC
    estimator in Figure~\ref{fig:IDoS_MoS2}. The step sizes used in the
    numerical differentiation are given in Table~\ref{tab:MoS2_MLMC}.
    The bottom right sub figure shows the density of states for all
    vacancy probabilities, $\Pvac$, together with that of the
    unperturbed material in an interval containing the bandgap of the
    unperturbed material.}
  \label{fig:DoS_MoS2}
\end{figure}

\section{Conclusions and future work}
\label{sec:conclusions}

We have studied Monte Carlo and MLMC sampling
methods for quantities of interest depending on the band structure of
quasi 2D materials with random impurities. 
We have presented a method of constructing control variates for the
quantities of interest by subdividing super cells into parts and using
the arithmetic mean of the quantity on the periodically extended
parts.
Using two tight-binding models on a honeycomb
lattice, we have empirically estimated the convergence rates of the
finite super cell bias, the variance on a finite super cell, and the
variance of the difference between a finite super cell sample and its
control variate, and found that for these test cases an MLMC approach
will be more computationally efficient than a standard 
Monte Carlo approach, which is in turn more efficient than using one
single sample on a larger super cell.

In the graphene test problem with a 32-by-32 super cell, a 2-level
Monte Carlo estimator of the same variance as a standard Monte Carlo
estimator was obtained at $61\%$ of the computational time of the
latter. This ratio should improve for a true Multilevel Monte Carlo
estimator as the size of the super cell increases. 
Indeed, in an $\MoS_2$ test problem, an MLMC estimator with
five super cell sizes ending with a 32-by-32 super cell, the estimated
computational savings were at least one order of magnitude.
More precisely,
based on the estimated convergence rates and costs, and on the asymptotic
complexity estimates, the work of an MLMC estimator
to meet accuracy $\tol$ in the quantity of interest in the test
problem is asymptotically proportional to $\tol^{-8/3}$ while 
the work of a standard Monte Carlo estimator with the same accuracy
grows like $\tol^{-10/3}$ as $\tol\to 0$.

Future work includes applying the MLMC approach for more demanding
quantities of interest, such as the conductivity tensor, other
geometries such as nano ribbons and bilayer heterostructures, 
studying more realistic distribution of
vacancies in the tight-binding model of $\MoS_2$, as well as
taking deformation of the lattice into account and using more accurate
density functional theory computations.

%
\bibliographystyle{siam}
\bibliography{report}

\end{document}